\documentclass[final]{siamltex}

\usepackage{amsmath}
\usepackage{graphicx}
\usepackage{latexsym}
\usepackage{amssymb}
\usepackage{amsbsy}
\usepackage{dsfont}
\DeclareSymbolFont{msbm}{U}{msb}{m}{n}
\DeclareMathSymbol{\C}{\mathalpha}{msbm}{'103}
\DeclareMathSymbol{\R}{\mathalpha}{msbm}{'122}
\DeclareMathSymbol{\Z}{\mathalpha}{msbm}{'132}
\DeclareMathSymbol{\N}{\mathalpha}{msbm}{'116}

\newcommand{\F}{\mathcal{F}}

\newcommand{\tA}{\tilde{A}}
\newcommand{\ta}{\tilde{a}}
\newcommand{\tw}{\tilde{w}}
\newcommand{\tc}{\tilde{c}}

\newtheorem{remark}{Remark}

\def\F{{\cal F}}
\def\U{{\cal U}}
\def\RR{\mathbb R}

\def\be{\begin{equation}}
\def\ee{\end{equation}}
\def\bea{\begin{eqnarray}}
\def\ba{\begin{array}{l}\displaystyle}
\def\eea{\end{eqnarray}}
\def\ea{\end{array}}

\def\proof{{\bf Proof \par}}

\begin{document}
\title{Asymptotic preserving Implicit-Explicit Runge-Kutta methods for non linear kinetic equations}

\author{Giacomo Dimarco\thanks{Universit\'{e} de Toulouse; UPS, INSA, UT1, UTM;
CNRS, UMR 5219; Institut de Math\'{e}matiques de Toulouse; F-31062
Toulouse, France. ({\tt giacomo.dimarco@math.univ-toulouse.fr}).}
\and Lorenzo Pareschi\thanks{Mathematics Department, University of
Ferrara and CMCS, Ferrara, Italy ({\tt lorenzo.pareschi@unife.it}).}
}
\maketitle

\begin{abstract}
We discuss Implicit-Explicit (IMEX) Runge Kutta methods which are
particularly adapted to stiff kinetic equations of Boltzmann type.
We consider both the case of easy invertible collision operators and
the challenging case of Boltzmann collision operators. We give
sufficient conditions in order that such methods are asymptotic
preserving and asymptotically accurate. Their monotonicity
properties are also studied. In the case of the Boltzmann operator,
the methods are based on the introduction of a penalization
technique for the collision integral. This reformulation of the
collision operator permits to construct penalized IMEX schemes which
work uniformly for a wide range of relaxation times avoiding the
expensive implicit resolution of the collision operator. Finally we
show some numerical results which confirm the theoretical analysis.
\end{abstract}

\maketitle



\section{Introduction}
The numerical solution of Boltzmann-type equations close to fluid
regimes represents a real challenge for numerical methods. In these
regimes, in fact, the intermolecular collision rate grows
exponentially and the collisional time becomes very small. On the
other hand, the actual time scale for evolution is the fluid dynamic
time scale, which can be much larger than the collisional time. A
non dimensional measure of the importance of collision is given by
the Knudsen number which is large in the rarefied regions and small in
the fluid ones. Standard computational approaches lose their
efficiency due to the necessity of using very small time steps in
deterministic schemes or, equivalently, a large number of collisions
in probabilistic approaches. Unfortunately the use of implicit
solvers originates a prohibitive computational cost due to the high
dimensionality and the nonlinearity of the collision operator.

Several authors have tackled the above problem for the Boltzmann
equation in the recent past (see ~\cite{BLM, Filbet, dimarco6,
toscani} and the references therein). In summary, the possible
approaches that permit to overcome such a difficulty can be
subdivided into two main classes. Domain decomposition strategies
and asymptotic preserving schemes. The first class of methods
permits to avoid the problem of very small Knudsen number by
identifying the regions where it is possible to use the reduced
fluid model and the regions where the full kinetic model must be
solved. The literature is this direction has a long history we
recall here references~\cite{bourgat, KT} and some recent works by
Degond and coauthors~\cite{degond1, dimarco5}. A closely related
research approach combines stochastic and deterministic solvers in the different
regions by originating hybrid methods~\cite{dimarco2, dimarco3}.

Concerning the
asymptotic preserving strategies, these techniques permit to solve
the full problem in the entire domain for all choices of time steps
and Knudsen numbers. Along this direction we quote the pioneering papers by Coron and Perthame~\cite{CP} for the BGK model and by Gabetta, Pareschi and Toscani~\cite{toscani} where an explicit exponential technique for the full Boltzmann equation capable to avoid the costly inversion of the collision term has been developed. More recently several improvements to the above approaches has been presented, we refer to \cite{BLM, Filbet, dimarco6}.

In the present paper we develop Implicit-Explicit
Runge-Kutta methods~\cite{Ascher, CK, PRimex} which are
particularly efficient for stiff non linear kinetic equations. In particular we generalize the approach recently introduced in
\cite{Filbet, dimarco7}. First we consider the case where the implicit inversion of the collision term does not represent a problem, like for example the case of simple BGK operators. Asymptotic preservation properties and monotonicity are carefully studied and analyzed. Subsequently we deal with the challenging case of the full Boltzmann equation. To this aim, following~\cite{Filbet} we introduce a penalization strategy based on a decomposition of the gain
term of the collision operator into an equilibrium and a non
equilibrium part. This permits to derive new penalized IMEX schemes which keep the good asymptotic preservation properties of standard IMEX schemes by avoiding the costly inversion of the collision term.
Similarly to~\cite{PRimex, dimarco6}, we derive sufficient
conditions for asymptotic preservation and asymptotic accuracy.
Monotonicity property are also considered and discussed.

We emphasize that the penalized IMEX schemes here developed apply to any large system of stiff
ordinary differential equations of the form
\begin{equation}
Y'=F(Y)+\frac1{\varepsilon} R(Y),\quad Y(t_0)=Y_0,
\end{equation}
where $\varepsilon>0$ is a small parameter, $Y, F(Y)\in\R^N$ and the
non-linear operator $R(Y)$ is a dissipative relaxation operator~\cite{CLL}. Such operator is endowed with a $n\times N$ matrix
$Q$ of rank $n<N$ such that $QR(Y)=0$, $\forall\,\, Y$. This gives
a vector of $n$ conserved quantities $y=QY$. Solutions which
belong to the kernel of the operator $R(Y)=0$ are uniquely
determined by the conserved quantities $Y=E(y)$ and characterize
the manifold of local equilibria.

The methods here studied are based on the following decomposition
\begin{equation}
R(Y)=N(Y)+L(Y),
\end{equation}
where $N(Y)$ represents the non-dissipative non-linear part and $L(Y)$ is a
linear term such that $L(Y)=0$ implies $Y=E(y)$. For example $L(Y)=A(E(y)-Y)$ where $A>0$
is an estimate of the Jacobian of $R$ evaluated at equilibrium. Note
that, at variance with standard linearization techniques which
operate on the short time scale, the operator is linearized on the
asymptotically large time scale. This decomposition permits to apply IMEX techniques which
are implicit in the linear part and explicit in the non-linear
part. The use of such techniques, as we will see, permits to achieve
unconditionally stable and asymptotic preserving methods at the cost of an explicit scheme.


The rest of the paper is organized as follows. First in Section 2 we
recall some basic aspects on kinetic equations and their fluid dynamic limits.
The notion of asymptotic preservation is also introduced. In
Section 3, we consider standard IMEX Runge-Kutta
schemes applied to kinetic equations and derive conditions for  asymptotic preservation and
asymptotic accuracy. Next, in section 4 we
introduce the penalized IMEX schemes for the full Boltzmann model. We analyze their asymptotic preservation and asymptotic
accuracy properties. Monotonicity in the homogeneous case is also discussed. Finally some numerical examples of schemes up to third order and conclusions are reported in Sections 5 and 6.

\section{The Boltzmann equation and related kinetic equations}

We consider kinetic equations of the form~\cite{cercignani}
\be
\partial_t f + v\cdot\nabla_{x}f
=Q(f),\label{eq:1}\ee
with initial data $f(x,v,t)|_{t=0}=f_{0}(x,v)$. Here $f(x,v,t)$ is a non negative function
describing the time evolution of the distribution of particles with
velocity $v \in \R^{3}$ and position $x \in \Omega \subset \R^{d_x}$
at time $ t
> 0$. For notation simplicity in the sequel we will omit the
dependence of $f$ from the independent variables $x,v,t$ unless
strictly necessary. The operator $Q(f)$ characterizes the particles
interactions and in the case of the quadratic Boltzmann collision operator of rarefied gas dynamics it reads
 \be Q_{B}(f)=\int_{\RR^3\times S^2} B(|v-v_*|,n)
[f(v')f(v'_*)-f(v)f(v_*)]\,dv_*\,dn \label{eq:Q} \ee where \be
v'=v+\frac12(v-v_*)+\frac12|v-v_*|n,\quad
v'_*=v+\frac12(v-v_*)-\frac12|v-v_*|n, \ee and $B(|v-v_*|,n)$ is a
nonnegative collision kernel characterizing the details of the
collision. It is described by the following equation \[
B(|v-v_*|,n)=\sigma\left(\frac{(v-v_*)}{|v-v_*|}\cdot n\right)
|v-v_*|^\gamma,
\] with $\gamma \in [0,3)$. The case $\gamma=1$ is
referred to as hard spheres case, whereas the simplified situation
$\gamma=0$, is referred to as Maxwell case. Note that in most
applications the angle dependence is ignored and $\sigma$ is assumed
constant.

The operator $Q(f)$ is such that the local
conservation properties are satisfied \be\int_{\R^{3}} \phi(v) Q(f)\,
dv=:\langle \phi Q(f)\rangle=0 \label{eq:QC}\ee where
$\phi(v)=\left(1,v,\frac{|v|^2}{2}\right)^T$ are the collision invariants. In
addition it satisfies the entropy inequality \be
\frac{d}{dt}H(f) = \int_{\R^{3}} Q(f)\log f dv
\leq 0,\qquad H(f)=\int_{\R^{3}}f\log f\,dv. \label{eq:entropy} \ee
The functions such that
$Q(f)=0$ are the local Maxwellian equilibrium functions
\be M[f]=M(\rho,u,T)=\frac{\rho}{(2\pi
T)^{3/2}}\exp\left(\frac{-|u-v|^{2}}{2T}\right), \label{eq:M}\ee
where $\rho$, $u$, $T$ are the density, mean velocity and
temperature of the gas in the $x$-position and at time $t$ defined as
\be (\rho,\rho u,E)^T=\langle \phi f \rangle, \qquad
T=\frac1{3\rho}(E-\rho|u|^2). \ee Due to its computational complexity, the Boltzmann collision operator
$Q_{B}(f)$ is often replaced in applications by simpler operators, like the BGK operator which substitutes the binary interactions with a
relaxation towards the equilibrium of the form~\cite{BGK} \be
Q_{BGK}(f)=\mu(M[f]-f),\ee
where $\mu > 0$. The validity of this operator in
describing the physics of non equilibrium phenomena has been the
subject of many papers in the past~\cite{cercignani}.

\subsection{Fluid-limit and asymptotic-preserving methods}
If we rescale the space
and time variables in (\ref{eq:1}) as \be x'=\varepsilon x, \ \
t'=\varepsilon t,\ee and omit the primes to keep notation simple, we obtain \be
\partial_t f + v\cdot\nabla_{x}f
=\frac1{\varepsilon}Q(f)\label{eq:1b}\ee where $\varepsilon$ is the Knudsen number a non dimensional quantity directly proportional to the mean free pat between particles.

Now integrating (\ref{eq:1b}) against the collision
invariants in the velocity space leads to the following set of non
closed conservations laws \be
\partial_t \langle \phi f\rangle+\nabla_x
\langle v\phi f\rangle=0.\label{eq:macr}\ee

Close to fluid regimes, the mean free path between two collisions is
very small. In this situation, passing to the limit $\varepsilon\rightarrow 0$ we formally obtain $Q(f,f)=0$ from (\ref{eq:1b}) and so $f=M[f]$. Thus, at least formally, we
recover the closed hyperbolic system of compressible Euler equations
\be
\partial_t U+\nabla_x\cdot \F(U)=0
\label{eq:Euler} \ee with
\[
U=\langle \phi M[f]\rangle = (\rho,\rho u,E)^T,
\]
\[
\F(U)=\langle v\phi M[f]\rangle=(\rho u, \varrho u \otimes u+pI,
Eu+pu)^T,\quad p=\rho T,
\]
where $I$ is the identity matrix. Note that the above conclusions are
independent on the particular choice of $Q(f,f)$ provided it satisfies (\ref{eq:QC}) and admits Maxwellian of the form (\ref{eq:M}) as local equilibrium functions.

For small but non zero values of the Knudsen number, the evolution
equation for the moments can be derived by the so-called
Chapman-Enskog expansion or Hilbert expansion~\cite{cercignani}.
These approaches originate the compressible Navier-Stokes equations
as a second order approximation with respect to $\varepsilon$ to the
solution of the Boltzmann equation. In this case, however, the
particular choice of the collision operator influences the structure
of the limiting Navier-Stokes system.

\begin{figure}
\begin{center}
\setlength{\unitlength}{1.1cm}
\begin{picture}(4,5)(1,1)
\put(1,5){\large{$P^{\varepsilon}$}}
\put(1,2){\large{${P^{\varepsilon}_{\Delta t}}$}}
\put(0,3.5){$\Delta t\to 0$} \put(4.5,3.5){$\Delta t\to 0$}
\put(2.2,5.5){$\varepsilon\to 0$}
\put(2.2,1.5){${\varepsilon\to 0}$}
\put(1.2,2.5){\vector(0,1){2.2}} \put(1.6,2.1){\vector(1,0){2.2}}
\put(4,5){\large{$P^{0}$}} \put(4,2){\large{${P^{0}_{\Delta
t}}$}} \put(4.2,2.5){\vector(0,1){2.2}}
\put(1.6,5.1){\vector(1,0){2.2}}
\end{picture}
\end{center}
\caption{$P^{\varepsilon}$ is the original
singular perturbation problem and $P^{\varepsilon}_{\Delta t}$ its
numerical approximation characterized by a discretization
parameter $\Delta t$. The \emph{asymptotic-preserving (AP) property} corresponds to the request
that $P^{\varepsilon}_{\Delta t}$ is a consistent discretization of
$P^0$ as $\varepsilon\to 0$ independently of $\Delta
t$.
}
\label{fg:AP}
\end{figure}
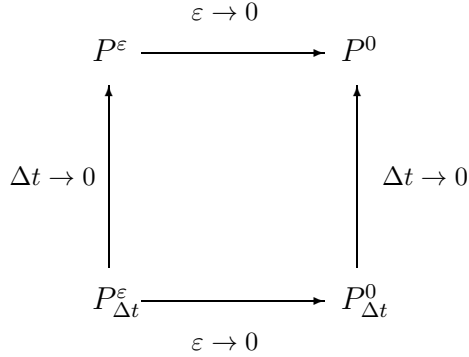

The construction of numerical schemes which are capable to capture the fluid-limit just described is closely connected with the notion of asymptotic-preserving schemes (see Figure \ref{fg:AP}). Here, in agreement with~\cite{Jin2, PRimex} we give the following definition of
asymptotic preserving methods for equation (\ref{eq:1})
\begin{definition}
A consistent time discretization method for (\ref{eq:1}) of stepsize
$\Delta t$ is {\em asymptotic preserving (AP)} if, independently of
the initial data and of the stepsize $\Delta t$, in
the limit $\varepsilon\to 0$ becomes a consistent time
discretization method for the reduced system (\ref{eq:Euler}).
\end{definition}

This definition does not imply that the scheme preserves the order
of accuracy in time in the stiff limit $\varepsilon\to 0$. In the
latter case we will say that the scheme is \emph{asymptotically accurate (AA)}.

In the sequel we will consider the development of asymptotic
preserving and asymptotically accurate schemes using the general
setting of IMEX Runge-Kutta methods.


\section{IMEX schemes for kinetic equations}

First we introduce the general formulation of IMEX schemes for kinetic equations together with some preliminary definitions.


An IMEX schemes applied to a kinetic equation of the form \be
\partial_t
f+v\cdot \nabla_x f=\frac1{\varepsilon} Q(f) \label{eq:BGK1}
\ee reads~\cite{PRimex}
\begin{eqnarray}
   F^{(i)} &=& \displaystyle f^{n}-\Delta t \sum_{j=1}^{i-1} \ta_{ij} v\cdot\nabla_x F^{(j)}+\Delta t\sum_{j=1}^{\nu} a_{ij}\frac{1}{\varepsilon} Q(F^{(j)})\label{eq:GIMEX} \\
   f^{n+1} &=& \displaystyle f^{n}-\Delta t
\sum_{i=1}^{\nu}\tw_{i}v\cdot\nabla_x F^{(i)}+\Delta
t\sum_{i=1}^{\nu}w_{i}\frac{1}{\varepsilon}Q(F^{(i)}).
\label{eq:GIMEX1}
\end{eqnarray}
The matrices $ \tA=(\ta_{ij} )$,
$\ta_{ij} = 0$ for $j \geq i$ and $A = (a_{ij})$ are
$\nu\times \nu$ matrices such that the resulting scheme is explicit
in $v\cdot\nabla_x f$, and implicit in $Q(f)$. In general, an IMEX
Runge-Kutta scheme, is characterized by the above defined two
matrices and the coefficient vectors $\tw =( \widetilde{
w}_{1},..,\tw_{\nu})^{T}$, $w =(w_{1},..,w_{\nu})^{T}$.
Since computational efficiency in the case of kinetic equations is of paramount importance, in the sequel we restrict our analysis to
diagonally implicit Runge-Kutta (DIRK) schemes for the source terms
($a_{ij} = 0,$ for $j
> i$). In fact, the use of a DIRK
scheme is enough to assure that the transport term $v\cdot\nabla_x f$ is always
evaluated explicitly. The type of schemes introduced can be
represented with a compact
notation by a double Butcher tableau 
\begin{table}[h!]
\begin{center}
\begin{tabular}{l|r c}
$\tc$ & $\tA$   \\
\hline\\
 & $\tw^{T}$
\end{tabular} \ \ \ \ \ \ \qquad
\begin{tabular}{l|r c}
$c$  & $A$   \\
\hline\\
 & $w^{T}$
\end{tabular}
\end{center}
\end{table}

\noindent
where the coefficients $\tc$ and $c$ are given by the
usual relation \be
\tc_{i}=\sum_{j=1}^{i-1}\ta_{ij}, \qquad
c_{i}=\sum_{j=1}^{i}a_{ij}.\ee
Using vector notations the schemes can be
written in compact form
\begin{eqnarray}
  F &=& \displaystyle f^{n}e+\Delta t \tA\,L(F)+\frac{\Delta t}{\varepsilon} A Q(F)\label{eq:GIMEXv} \\
  f^{n+1} &=& \displaystyle f^{n}+\Delta t
\tw^{T}L(F)+\frac{\Delta
t}{\varepsilon}w^{T}Q(F),
\label{eq:GIMEXv1}
\end{eqnarray}
where $e=(1,1,..,1)^{T}\in \R^{\nu}$, $F=(F^{(1)},\ldots,F^{(\nu)})^T$, $Q(F)=(Q(F^{(1)}),\ldots,Q(F^{(\nu)}))^T$ and $L(F)=(L(F^{(1)}),\ldots,L(F^{(\nu)}))$ with $L(F^{(i)})=-v\cdot\nabla_{x}F^{(i)}$.

We refer to~\cite{Ascher, CK, HairerWanner, PRimex} for more details on
the order conditions for IMEX schemes. Let us remark that IMEX schemes are a particular case of additive Runge-Kutta methods and so the order conditions can be derived as a generalization of the notion of Butcher tree~\cite{HairerWanner}. In particular, under the assumptions $\widetilde c=c$ and
$\widetilde w=w$, mixed order conditions are automatically satisfied up to third order.

Before stating the main properties concerning asymptotic preservation it is useful to characterize the different IMEX schemes we will consider in the sequel accordingly to the structure of the DIRK method. Following~\cite{BPR} we have
\begin{definition} We call an IMEX-RK method of \emph{type A} (see \cite{PRimex}) if the matrix $A \in \RR^{\nu \times \nu}$ is invertible, or equivalently $a_{ii}\neq 0$, $i=1,\ldots,\nu$.
We call an IMEX-RK method of \emph{type CK} (see \cite{CK}) if the matrix $A$ can
be written as
\be
A = \left(\begin{array}{ll} 0 & 0\\
                        a  &  \hat{A}\end{array}\right),
\label{CK1}
\ee
with $a=(a_{21},\ldots,a_{\nu 1})^T\in \RR^{(\nu-1)}$ and the submatrix $\hat{A} \in \RR^{(\nu-1)
\ \times \ (\nu-1)}$ invertible, or equivalently $a_{ii}\neq 0$, $i=2,\ldots,\nu$. In the special case $a=0$, $w_1=0$ the
scheme is said to be of \emph{type ARS} (see \cite{Ascher}) and the DIRK method is reducible to a method using $\nu-1$ stages.
\end{definition}

We will also make use of the following representation of the matrix $\tA$ in the explicit Runge-Kutta method
\be
\tA=
  \left(
          \begin{array}{ll}
            0 & 0  \\
            \ta & \hat{\tA}\end{array}
        \right),
\label{CK2}\ee
where $\ta=(\ta_{21},\ldots,\ta_{\nu 1})^T\in\R^{\nu-1}$ and $\hat{\tA}\in\R^{\nu-1\times\nu-1}$.

The following definition will be also useful to characterize the properties of the methods in the sequel.
\begin{definition}
We call an IMEX-RK method \emph{implicitly stiffly accurate (ISA)} if the corresponding DIRK method is \emph{stiffly accurate}, namely
\be
a_{\nu i}=w_i,\quad i=1,\ldots,\nu.
\ee
If in addition the explicit methods satisfies
\be
\ta_{\nu i}=\tw_i,\quad i=1,\ldots,\nu
\ee
the IMEX-RK method is said to be \emph{globally stiffly accurate (GSA)} or simply \emph{stiffly accurate}.
\end{definition}

Some remarks are in order.
\begin{remark}
\begin{itemize}
\item For type $A$ IMEX schemes we have $a_{11}\neq 0$ and $\ta_{11}=0$ thus $c_1\neq \tc_1$ and we cannot assume the simplifying condition $c=\tc$.
\item For IMEX schemes of type $A$ or type $CK$ the GSA property implies $\tw_{\nu}=0$ and $w_{\nu}\neq 0$ thus we cannot assume the simplifying condition $w=\tw$. Note that for GSA schemes the numerical solution is the same as the last stage value, namely $f^{n+1}=F^{(\nu)}$.
\item From the observations above it is clear that order conditions for GSA type $A$ IMEX schemes are particularly restrictive since $\tc\neq c$ and $\tw\neq w$.
\end{itemize}
\end{remark}

\subsection{Asymptotic preserving IMEX schemes}

We give now conditions for an IMEX scheme to satisfy asymptotic
preservation and asymptotic accuracy. Here we do not consider the computational challenges related to the inversion of the implicit collision operator $Q(f)$. We will focus on these aspects in the second part of the paper.

We can state the following theorem which show that type $A$ IMEX schemes are asymptotic preserving and
asymptotically accurate.
\begin{theorem}
\label{th:ap} If the IMEX method is of type $A$ then in the limit
$\varepsilon\rightarrow 0$, scheme
(\ref{eq:GIMEXv})-(\ref{eq:GIMEXv1}) becomes the explicit RK scheme characterized by
($\tA, \tw, \tc$) applied to the limit
Euler system (\ref{eq:Euler}).
\end{theorem}
\newpage \proof To prove the above theorem let us first multiply the IMEX
method (\ref{eq:GIMEXv})-(\ref{eq:GIMEXv1}) by the collision
invariants $\phi(v)=1, v, v^2$ and integrate the result in velocity
space. We obtain the explicit Runge-Kutta methods applied to the
moment system (\ref{eq:macr})
\begin{eqnarray}
\displaystyle   \langle \phi F\rangle &=& \langle \phi f^{n}e\rangle+\Delta t \tA\langle \phi L(F)\rangle
\label{eq:GIMEXfm} \\
\displaystyle   \langle \phi f^{n+1}\rangle &=& \langle\phi f^{n}\rangle+\Delta t
\tw^{T}\langle\phi L(F)\rangle.
\label{eq:GIMEX1fm}
\end{eqnarray}
Now let us rewrite equation (\ref{eq:GIMEXv})
in the following form \be \varepsilon F=\varepsilon Fe
+\varepsilon\Delta t \tA L(F)+\Delta t A Q(F).\ee
Since $A$ is invertible we can solve for $Q(F)$ to get \be \Delta t Q(F)=\varepsilon
A^{-1}\left(F-f^{n}e-\Delta t
\tA L(F)\right).\label{eq:th1}\ee
As
$\varepsilon\rightarrow 0$ we get \be \Delta t Q(F)=0 \Rightarrow F=M[F].\ee Thus (\ref{eq:GIMEXfm})-(\ref{eq:GIMEX1fm})
 becomes the explicit Runge-Kutta method applied to the limiting Euler system (\ref{eq:Euler})
\begin{eqnarray}
\displaystyle   {\cal U} &=& U^{n}e-\Delta t \tA\nabla_x\cdot \F({\cal U})
\label{eq:GIMEXfeuler} \\
\displaystyle   U^{n+1} &=& U^{n}+\Delta t
\tw^{T}\nabla_x\cdot \F({\cal U}),
\label{eq:GIMEX1feuler}
\end{eqnarray}
where ${\cal U}=(\U^{(1)},\ldots,\U^{(\nu)})^T$, $\F({\cal U})=(\F(\U^{(1)}),\ldots,\F(\U^{(\nu)}))^T$, $\U^{(i)}=\langle \phi M[F^{(i)}]\rangle$ and $\F(\U^{(i)})=\langle \phi L(M[F^{(i)}])\rangle$.

To conclude the proof we must be able to pass to the limit $\varepsilon\to 0$ in the numerical solution. This is in fact possible since for type $A$ IMEX schemes the numerical solution is
independent on $\varepsilon$. In fact, we have \be
f^{n+1}=f^{n}+\Delta t\tw^{T}L(F)+\frac{\Delta
t{w}^{T}}{\varepsilon}Q(F)\ee which thanks to (\ref{eq:th1}) yields \[
f^{n+1}=f^{n}+\Delta
t\tw^{T}L(F)+{w}^{T}A^{-1}\left(F-f^{n}e-\Delta
t\tA L(F)\right)\] and finally \be
f^{n+1}=f^{n}\left(1-w^{T}A^{-1}e\right)+\Delta
t\left(\tw^{T}-w^{T}A^{-1}\tA\right)L(F)+{w}^{T}A^{-1}F.
\label{eq:nsA}\ee\\
$\Box$

An important property of the schemes is that in
the limit $\varepsilon\rightarrow 0$ the distribution function is
projected over the equilibrium $f^{n+1}\rightarrow M[f^{n+1}]$. From (\ref{eq:nsA}) it is clear that this property is achieved if the following conditions are satisfied
\be
            w^{T}A^{-1}e=1,\qquad
            \tw^{T}=w^{T}A^{-1}\tA, \qquad
            w^{T}A^{-1}M[F]=M[f^{n+1}].
        \label{eq:fn1}
\ee
The first condition corresponds to the classical $L$-stability requirement of the DIRK method~\cite{HairerWanner}.
Since the third condition depends on the stage values vector $F$ the only possibility to satisfy (\ref{eq:fn1}) is that the IMEX scheme is GSA. In this case, in fact, we have
\[
w^T A^{-1}=(0,\ldots,0,1)^T,\quad M[F^{\nu}]=M[f^{n+1}].
\]
Thus, we can state
\begin{theorem}
If the IMEX scheme is of type A and GSA then
\begin{equation}
\lim_{\varepsilon\to 0} f^{n+1}=M[f^{n+1}].
\end{equation}
\end{theorem}
As observed, the request that the
matrix $A$ is invertible can be quite restrictive for high order methods since we cannot assume the simplifying condition $c=\tc$. However, under additional hypothesis, we can obtain schemes
which are asymptotic preserving and asymptotically accurate even when the first row of $A$ contains only zeros.

In order to do this we first introduce the notion of initial
data {consistent} with the limit problem.
\begin{definition} The initial data for  equation
(\ref{eq:BGK1}) are said
\emph{consistent} or \emph{well prepared} if
\be f_0(x,v)=M[f_0(x,v)]+g^{\varepsilon}(x,v),\qquad \lim_{\varepsilon\to 0} g^{\varepsilon}(x,v)=0.\label{eq:cid}\ee
\end{definition}
We can now state the following
\begin{theorem}
\label{th:apCK} If the IMEX scheme is of type CK and GSA then for consistent initial
data
 in the limit $\varepsilon\rightarrow 0$ scheme
(\ref{eq:GIMEXv})-(\ref{eq:GIMEXv1}) becomes the explicit RK scheme characterized by
($\tA, \tw, \tc$) applied to the limit
Euler system (\ref{eq:Euler}). \end{theorem}

\proof
First we rewrite $F=\left(F^{(1)},
\hat{F}\right)^{T}$, $w=(w_1,\hat{w})^T$, $\tw=(\tw_1,\hat{\tw})^T$ with $\hat{F},\hat{w},\hat{\tw}\in\R^{\nu-1}$ and the matrix $\tA$ of the explicit tableau as in (\ref{CK2}).

Then for type CK the IMEX schemes (\ref{eq:GIMEXv})-(\ref{eq:GIMEXv1}) can be written as
\begin{eqnarray}
\nonumber
F^{(1)}&=&f^{n}\\[-.25cm]
\label{eq:theo3a}
\\[-.25cm]\nonumber
\hat{F} &=& \displaystyle f^{n}\hat{e}+\Delta t \ta L(F^{(1)})+\Delta t
\hat{\tA}L(\hat{F})+\frac{\Delta
t}{\varepsilon} a Q(F^{(1)})+
\frac{\Delta t}{\varepsilon}
\hat{A}Q(\hat{F}),\\
f^{n+1}&=&f^n +\Delta t\tw_1 L(F^{(1)})+\Delta t\hat{\tw}^T L(\hat{F})+w_1\frac{\Delta t}{\varepsilon} Q(F^{(1)})+\hat{w}^T\frac{\Delta t}{\varepsilon} Q(\hat{F})
\label{eq:theo3b}
\end{eqnarray}
where $\hat{e}=(1,..,1)\in \mathbb{R}^{\nu-1}$. The corresponding moment scheme yields
\begin{eqnarray}
\nonumber
\langle \phi F^{(1)}\rangle &=& \langle \phi f^{n}\rangle\\[-.3cm]
\label{eq:GIMEXaa} \\[-.3cm]\nonumber\displaystyle   \langle \phi \hat{F}\rangle &=& \langle \phi f^{n}\hat{e}\rangle+\Delta t \ta \langle\phi L(F^{(1)})\rangle
+\Delta t
\hat{\tA}\langle\phi L(\hat{F})\rangle\\
\displaystyle   \langle \phi f^{n+1}\rangle &=& \langle\phi f^{n}\rangle+\Delta t\tw_1\langle\phi L(F^{(1)})\rangle+\Delta t\hat{\tw}^T \langle\phi L(\hat{F})\rangle.
\label{eq:GIMEXbb}
\end{eqnarray}

Solving now,
the second equation in (\ref{eq:theo3a}) for $Q(\hat{F})$ we get \be
\Delta t Q(\hat{F})=\varepsilon
\hat{A}^{-1}\left[\hat{F}-f^{n}\hat{e}-\Delta t
\hat{\tA}L(\hat{F})-\Delta t \ta
L(f^{n})\right]-\Delta t
\hat{A}^{-1}a Q(f^{n}).
\label{eq:ll1}\ee
As $\varepsilon\to 0$ we obtain
\be
Q(\hat{F})= -
\hat{A}^{-1}a Q(f^{n}).
\label{eq:limit}
\ee
Since $f^n=M[f^n]+g^\varepsilon$ in the limit
$\varepsilon\rightarrow 0$ we have $f^n=M[f^n]$, $Q(f^n)=0$, and (\ref{eq:limit}) reduces to $Q(\hat{F})=0$ which implies $\hat{F}=M[\hat{F}]$.
Moreover, since the scheme is GSA, we also have $f^{n+1}=M[f^{n+1}]$ and thus at the next time step the initial value remains consistent and the moments system (\ref{eq:GIMEXaa})-(\ref{eq:GIMEXbb}) corresponds to the explicit Runge-Kutta methods for the Euler equations.\\
$\Box$

If we now remove the assumptions of GSA and consistent initial data, the numerical solution takes the form
\begin{eqnarray}
\nonumber
f^{n+1}&=&f^n\left(1-\hat{w}^T
\hat{A}^{-1}\hat{e}\right)+\Delta t\left(\tw_1-\hat{w}^T
\hat{A}^{-1}{\ta}\right) L(f^{n})+\Delta t\left(\hat{\tw}^T-\hat{w}^T
\hat{A}^{-1}\hat{\tA}\right) L(\hat{F})\\[-.25cm]
\label{eq:ns}
\\[-.25cm]
\nonumber
&+&\hat{w}^T
\hat{A}^{-1}\hat{F}+\frac{\Delta t}{\varepsilon}\left(w_1-\hat{w}^T \hat{A}^{-1}a\right)Q(f^n).
\end{eqnarray}
Thus in order to be able to pass to the limit $\varepsilon\to 0$ we require the additional condition
\be
w_1=\hat{w}^T \hat{A}^{-1}a.
\ee
The above equation is obviously satisfied by IMEX schemes of type ARS since $w_1=0$ and $a=0$. In addition, if the IMEX scheme is of type ARS equation (\ref{eq:limit}) reduces to $Q(\hat{F})=0$ which implies $\hat{F}=M[\hat{F}]$. Thus, except for the initial layer effect given by the $O(\Delta t)$ term $L(f^n)$ in (\ref{eq:GIMEXaa})-(\ref{eq:GIMEXbb}) IMEX scheme of type ARS, even if not necessarily asymptotically accurate, are asymptotic preserving without further assumptions. The additional requirement $\tw_1=0$ in this case suffices to give an asymptotically accurate method.

Finally, we can state the following results on the asymptotic behavior of the numerical solution
\begin{theorem}
If the IMEX scheme is of type CK and GSA then
\be
\lim_{\varepsilon\to 0} f^{n+1}=M[f^{n+1}],
\ee
if one of the following conditions is satisfied
\begin{itemize}
\item[\rm (a)] the initial
data is consistent;
\item[\rm (b)]
\be
\hat{e}_{\nu}^T\hat{A}^{-1}a=0,
\label{eq:app}
\ee
where $\hat{e}_{\nu}=(0,\ldots,0,1)^T\in\R^{\nu-1}$.
\end{itemize}
\end{theorem}
Note that (\ref{eq:app}) is true for ARS type schemes.
We end this paragraph with some remarks
\begin{remark}
\begin{itemize}
\item In case the IMEX methods are not GSA the numerical solution may originate a final layer effect which leads to reduction of accuracy in the kinetic variable $f^{n+1}$. Similarly for type CK IMEX schemes, in case the initial data is not consistent,
the solution may exhibit an initial layer which gives rise to a reduction of accuracy in the numerical solution.
These phenomena can be cured using smaller time steps or extrapolation techniques only for the very last, respectively first, step of the computation.
\item To avoid the loss of accuracy in the kinetic variable of non GSA IMEX methods it is also possible to impose the additional conditions up to third order~\cite{BPR}
\begin{eqnarray}\label{oc_index1}
b^T A^{-1} \tilde{c}  = 1, \ \ \ b^T A^{-1} \tilde{c}^2  = 1, \ \ \ b^T A^{-1} \tilde{A}\tilde{c}  = 1/2,
\end{eqnarray}
where $ \tilde{c} = \tilde{A}e$. For the
type CK (consequently for the type ARS) we can rewrite these
algebraic order conditions replacing $A^{-1}$ with $\hat{A}^{-1}$.
\end{itemize}
\end{remark}

\subsection{Properties of IMEX schemes for relaxation operators}
\label{par:BGK}
In this paragraph we consider the particular case of BGK relaxation operators of the form $Q_{BGK}(f)=\mu(M[f]-f)$.
A fundamental property of the IMEX scheme
(\ref{eq:GIMEXv})-(\ref{eq:GIMEXv1}) applied to relaxation operators is that it can be solved explicitly
despite the nonlinearity of $M[f]$.

To understand this let us remark that since the implicit scheme is a DIRK method equation (\ref{eq:GIMEX}) takes the form
\begin{eqnarray}
   F^{(i)} = \displaystyle f^{n}+\Delta t \sum_{j=1}^{i-1} \ta_{ij} L(F^{(j)})+\Delta t\sum_{j=1}^{i} a_{ij}\frac{\mu}{\varepsilon}(M[F^{(j)}]-F^{(j)})\label{eq:GIMEXdirk}
\end{eqnarray}
where the only implicit term is the diagonal factor
$M[F^{(i)}]-F^{(i)}$ in which $M[F^{(i)}]$ depends only on the
moments $\langle\phi F^{(i)}\rangle$. If we now integrate equation
(\ref{eq:GIMEXdirk}) against the collision invariants thanks to the
conservations (\ref{eq:QC}) we obtain
\begin{eqnarray}
 \langle \phi F^{(i)}\rangle = \displaystyle \langle \phi f^{n}\rangle +\Delta t \sum_{j=1}^{i-1} \ta_{ij} \phi L(F^{(j)})\rangle\label{eq:GIMEXdirkm}
\end{eqnarray}
which corresponds to the stages of the explicit Runge-Kutta method
applied to the moment system (\ref{eq:macr}). Thus $\langle \phi
F^{(i)}\rangle$, and so $M[F^{(i)}]$, can be explicitly evaluated
and system (\ref{eq:GIMEXdirk}) is explicitly solvable. This
property has been used, for example, in~\cite{Puppo} to implement
efficiently IMEX methods for the BGK system.

Monotonicity properties for Runge-Kutta methods have been studied by
several authors in the recent past~\cite{Fer, GST}. The same
properties have been analyzed in the specific case of additive Runge
Kutta schemes in~\cite{Hi1}. Here, we restrict our study to the
space homogeneous case ($L(f)=0$) leaving to a future research the
analysis in the non homogeneous case.  Even if this case is
particularly simple, we report briefly some results, since they are
the basis of the analysis for the Boltzmann equation that we will
discuss in the next section. Note, however, as shown in the above
cited papers, that monotonicity conditions in the non homogeneous
case are rather restrictive. Thus, if positivity is strictly
demanded, as for instance in the case of Monte Carlo methods,
splitting strategies are normally preferable~\cite{dimarco6}.


In the homogeneous case the method reduces to a simple application of the implicit scheme and reads
\begin{eqnarray}
\label{eq:s1}
F=f^{n}e+\frac{\mu\Delta t}{\varepsilon}A(M[f^n]e-F),\\
\label{eq:s2}
f^{n+1}=f^n+w^T\frac{\mu\Delta t}{\varepsilon}(M[f^n]e-F),
\end{eqnarray}
where now the local equilibrium $M[F]=M[f^n]e$ is independent of time since $\langle\phi f^{n+1}\rangle=\langle\phi f^n\rangle$.

Let us
define $z=\mu\Delta t/\varepsilon$ and solve for $F$, we get \begin{eqnarray}
F&=&(I+zA)^{-1}(f^{n}e+zAM[f^n]e)\\
f^{n+1}&=&(1-w^{T}A^{-1}e)f^{n}+w^{T}A^{-1}F,
\end{eqnarray}
where we used the fact that
\be
z(M[f^n]e-F)=A^{-1}[F-f^{n}e].\label{eq:simp}\ee
Thus we have the following
\begin{proposition}
Sufficient conditions to guarantee that $f^{n+1}\geq 0$ when $f^n\geq 0$ in (\ref{eq:s1})-(\ref{eq:s2}) are that
\begin{eqnarray}
\label{eq:posn}
\left(I+zA\right)^{-1}e\geq 0,&\qquad& \left(I+zA\right)^{-1}Ae\geq 0,\\
1-w^{T}A^{-1}e \geq 0,&\qquad& w^{T}A^{-1}\geq 0.
\label{eq:posn1}
\end{eqnarray}
\end{proposition}
Note that conditions (\ref{eq:posn})-(\ref{eq:posn1}) must be
interpreted component by component. In particular, (\ref{eq:posn})
depend on $z$ and (\ref{eq:posn1}) require $A$ to be invertible. For
ISA schemes conditions reduce to (\ref{eq:posn}). The maximum value
of $\bar{z}\geq 0$ such that conditions (\ref{eq:posn}) are
satisfied for $z\in [0,\bar{z}]$ is usually referred to as the
\emph{radius of absolute monotonicity} of the scheme~\cite{Fer}.

Let us remark that since
\be
(F-M[f^n]e)=(I+zA)^{-1}(f^{n}-M[f^n])e.\label{eq:simpp}\ee
the numerical solution can be written as
\be
f^{n+1}=R(z)f^{n}+(1-R(z))M[f^n],
\label{eq:dirks}
\ee
where $R(z)=1-zw^{T}(I+zA)^{-1}e$ is the stability function of the DIRK method~\cite{HairerWanner}. Note that (\ref{eq:dirks}) has the same structure of the exact solution to our problem where $R(z)$ is an approximation of $\exp(-z)$.
Then, we can state
\begin{proposition}
A sufficient condition to guarantee that $f^{n+1}\geq 0$ when $f^n\geq 0$ in (\ref{eq:s1})-(\ref{eq:s2}) is that $0\leq R(z) \leq 1$.
\end{proposition}

Clearly the above proposition defines a subset of the absolute
stability region of the method~\cite{HairerWanner} characterized by
$|R(z)|\leq 1$. Moreover, being based on a convexity argument, it
gives a sufficient condition for the entropic property of the
scheme, namely $H(f^{n+1})\leq H(f^n)$ (see (\ref{eq:entropy})).
\begin{remark} In this simple situation the L-stability of the implicit methods implies the AP property. In fact if $R(z)\to R(\infty)$ as $z\to\infty$, then $R(\infty)=0$ yields $f^{n+1}\to M[f^n]$ as $\varepsilon\to 0$. A  recursive computation also shows that
\be
f^{n+1}=R(z)^{n+1} f_0+(1-R(z)^{n+1})M[f_0],
\label{eq:Rrec}
\ee
since $M[f^n]=M[f_0]$, $\forall\, n$. Therefore $f^{n+1}\to M[f_0]$ as $n\to\infty$ provided $|R(z)|<1$. Hence if $0<|R(\infty)|<1$ we have the weaker AP property $f^{n+1}\to M[f_0]$ as $\varepsilon\to 0$, $n\to \infty$.
\end{remark}

\section{Penalized IMEX schemes for the Boltzmann
equation} In the sequel, we will describe how to modify the IMEX
Runge-Kutta strategy in the challenging case of the
Boltzmann equation (\ref{eq:1}) where the collisions are described
by the operator $Q(f)=Q_{B}(f)$ introduced in (\ref{eq:Q}).

Although the approach just described remains formally valid, in
practice a major difficulty concerns the need to solve the system of
nonlinear equations originated by the application of the DIRK method
to the collision operator $Q_B(f)$. The computational cost of such
integral operator, characterized by a five fold nonlinear integral
which depends on the seven dimensional space $(x,v,t)$, is extremely
high and makes it prohibitive the use of iterative solvers. Thus, to
overcome this difficulty, the idea is to reformulate the collision
part using a suitable penalization term. This idea, using a BGK
model as penalization term, has been introduced recently
in~\cite{Filbet}. We refer to~\cite{dimarco6} for an extension of
this approach to exponential Runge-Kutta techniques.

\subsection{Penalization of the collision integral}
The key idea of penalization techniques for the Boltzmann collision
operator in stiff regimes is to observe that in such regimes the
solution $f$ is close to the Maxwellian equilibrium $M[f]$. Due to
its computational and analytical complexity a large variety of
simplified approximations to the Boltzmann collision operator which
are valid close to equilibrium have been derived in the
literature~\cite{BGK, cercignani, Hol}. We mention here the BGK
approximation~\cite{BGK} $Q_{BGK}(f)=\mu(M[f]-f)$ already discussed,
which is capable to capture the leading order term, i.e. the
compressible Euler equations. We mention also the ES-BGK relaxation
approximation~\cite{Hol} $Q_{ES}(f)=\mu(M_{\varepsilon}[f]-f)$ where
a modified equilibrium $M_{\varepsilon}[f]$, $M_0[f]=M[f]$ is used,
such that the equation has the advantage of matching the
$O(\varepsilon)$ expansion, i.e. the compressible Navier-Stokes
system. Finally, another possibility is to consider the linear
Boltzmann equation in the form~\cite{cercignani} \be
Q_{LB}(f)=\int_{\RR^3\times S^2} B(|v-v_*|,n)
\{f(v')M[f](v'_*)-f(v)M[f](v_*)\}\,dv_*\,dn. \label{eq:QL} \ee A
fundamental property shared by all the previous operators is that
their kernel is spanned by the local Maxwellian equilibrium $M[f]$,
so that the $\varepsilon\to 0$ asymptotic behavior of the original
Boltzmann equation is preserved.

Let us now denote with $Q_{P}(f)$ a general operator which will be
used to penalize the original Boltzmann operator $Q_{B}(f)$. As
discussed, the characteristics of $Q_P(f)$ are to be computable and
invertible at a low computational cost and that it preserves the
local equilibrium, namely $Q_P(f)=0$ implies $f=M[f]$.

We will then rewrite the collision operator in the form
\be
Q_B(f)=(Q_B(f)-Q_P(f))+Q_P(f)=G_P(f)+Q_P(f),
\label{eq:QP}
\ee
where by construction $\langle\phi\, G_P(f)\rangle=0$, and the corresponding kinetic equation reads
\be
\partial_t
f+v\cdot\nabla_x
f=\frac{1}{\varepsilon}G_P(f)+\frac{1}{\varepsilon}Q_P(f).\label{eq:Bref2}
\ee
Now, the idea is to use a numerical scheme in which only
the simpler operator $Q_P(f)$ is treated implicitly, while the
term $G_P(f)$ describing the deviations of the true Boltzmann operator $Q_B(f)$ from the simplified operator $Q_P(f)$ and the
convection term $\nabla_x f$ are treated explicitly.  This
however, as we will see, introduces some additional stability
requirements in order for the IMEX schemes to preserve the
asymptotic behavior of the equation.

Let us remark that the above penalization method in the case
$Q_P(f)=Q_{BGK}(f)$ admits an interesting interpretation that
relates the previous approach with decomposition
methods~\cite{dimarco2,dimarco3,Lem}. First, we observe that we can
rewrite $Q_{B}(f)$ as \be Q_{B}(f)=P(f)-\mu f,\label{eq:10} \ee
where $P(f)=Q_{B}(f)+\mu f$ and $\mu>0$ is a constant such that
$P(f)\geq 0$. Typically $\mu$ is an estimate of the largest rate of
the negative term in the Boltzmann operator
 \be \mu\geq \int_{\RR^3\times S^2} B(|v-v_*|,n)
f(v_*)\,dv_*\,dn. \label{eq:P}\ee By construction, the following property is
verified by the operator $P(f,f)$\be \frac1{\mu}\langle \phi
P(f,f)\rangle = \langle \phi f\rangle=U. \ee Thus $P(f,f)/\mu$ is a
density function and we can consider the following decomposition \be
P(f,f)/\mu=M[f]+g,\label{eq:dec}\ee where the function $g=P(f,f)/\mu-M[f]$
represents the non equilibrium part of the distribution function and
from the definition above it follows that $g$ is in general non
positive. Moreover since $P(f,f)/\mu$ and $M[f]$ share the first
three moments we have $\langle \phi g\rangle=0$.

Thus the collision operator can be rewritten in the form \be
Q_{B}(f,f)=\mu M[f]+\mu g-\mu f = {\mu}g+{\mu}(M[f]-f)=G_{BGK}(f)+Q_{BGK}(f),\label{eq:11}\ee
which is the same as (\ref{eq:QP}) using the BGK model as penalizer.
As we will see, thanks to the derivation above, the use of the BGK model is particularly interesting to study the monotonicity properties of the penalized IMEX schemes.

\subsection{Asymptotic preserving penalized IMEX schemes} We can now
introduce the general class of penalized IMEX Runge-Kutta schemes for the
Boltzmann equation in the form
\begin{eqnarray}
&   F^{(i)} = \displaystyle f^{n}+\Delta t \sum_{j=1}^{i-1} \ta_{ij}
\left(\frac{1}{\varepsilon} G_P(F^{(j)})-v\cdot\nabla_x F^{(j)}\right)+\Delta t\sum_{j=1}^{\nu}a_{ij}\frac{1}{\varepsilon} Q_P(F^{(j)})\label{eq:GIMEXb} \\
&   f^{n+1} = \displaystyle f^{n}+\Delta t
\sum_{i=1}^{\nu}\tw_{i}\left(\frac{1}{\varepsilon}
G_P(F^{(i)})-v\cdot\nabla_x F^{(i)}\right)+\Delta
t\sum_{j=1}^{\nu}w_i\frac{1}{\varepsilon}Q_P(F^{(i)}),
\label{eq:GIMEX1b}
\end{eqnarray}
or equivalently, using vector notation as
\begin{eqnarray}
&   F = \displaystyle f^{n}e+\Delta t \tA\left(\frac{1}{\varepsilon} G_P(F)
+L(F)\right)+\Delta t A\frac{1}{\varepsilon}Q_P(F)\label{eq:GIMEXbv} \\
&   f^{n+1} = \displaystyle f^{n}+\Delta t
\tw^{T}\left(\frac{1}{\varepsilon}
G_P(F)+L(F)\right)+\Delta t w^{T}\frac{1}{\varepsilon}Q_P(F),
\label{eq:GIMEX1bv}
\end{eqnarray}
where $G_P(F)=(G_P(F^{(1)}),\ldots,G_P(F^{(\nu)}))^T$.

We want now to derive conditions for the penalized IMEX schemes to be AP and
asymptotically accurate. As we will see, even for type A IMEX schemes, additional conditions are required to achieve asymptotic preservation.
\begin{theorem}
\label{th:pap} If the penalized IMEX method is of type $A$ and satisfies
\be
\tw^{T}=w^{T}A^{-1}\tA,
\label{eq:ai1}
\ee
 then in the limit
$\varepsilon\rightarrow 0$, scheme
(\ref{eq:GIMEXbv})-(\ref{eq:GIMEX1bv}) becomes the explicit RK scheme characterized by
($\tA, \tw, \tc$) applied to the limit
Euler system (\ref{eq:Euler}).
\end{theorem}

\proof
First let us observe that if we multiply the penalized IMEX method
(\ref{eq:GIMEXbv})-(\ref{eq:GIMEX1bv}) by the collision invariants
$\phi(v)=1, v, v^2$ and integrate the result in velocity space we obtain again the explicit Runge-Kutta methods applied to the moment system (\ref{eq:macr})
\begin{eqnarray}
\displaystyle   \langle \phi F\rangle &=& \langle \phi f^{n}e\rangle+\Delta t \tA\langle \phi L(F)\rangle
\label{eq:GIMEXfmb} \\
\displaystyle   \langle \phi f^{n+1}\rangle &=& \langle\phi f^{n}\rangle+\Delta t
\tw^{T}\langle\phi L(F)\rangle.
\label{eq:GIMEX1fmb}
\end{eqnarray}

Now let us rewrite equation (\ref{eq:GIMEXv})
in the following form \be \varepsilon F=\varepsilon Fe
+\varepsilon\Delta t \tA \left(\frac{1}{\varepsilon} G_P(F)
+L(F)\right)+\Delta t A Q_P(F).\ee
Since $A$ is invertible we can solve for $Q_P(F)$ to get \be \Delta t Q_P(F)=\varepsilon
A^{-1}\left(F-f^{n}e-\Delta t
\tA \left(\frac{1}{\varepsilon} G_P(F)
+L(F)\right)\right).\label{eq:th1b}
\ee
As
$\varepsilon\rightarrow 0$ we get \be Q_P(F)=-A^{-1}\tA G_P(F).\ee
Since $A^{-1}\tA$ is lower triangular with diagonal elements equal to zero we get
\be
Q_P(F^{(i)})=0\quad\Rightarrow\quad F^{(i)}=M[F^{(i)}],\quad i=1,\ldots,\nu.
\ee
Thus (\ref{eq:GIMEXfmb})-(\ref{eq:GIMEX1fmb})
 becomes the explicit Runge-Kutta method applied to the limiting Euler system (\ref{eq:Euler}).

Unfortunately, now the numerical solution still depends on $\varepsilon$. In fact, inserting (\ref{eq:th1b}) into (\ref{eq:GIMEX1bv}) we obtain
\[
f^{n+1}=f^{n}+\Delta
t\tw^{T}\left(\frac{1}{\varepsilon}
G_P(F)+L(F)\right)+{w}^{T}A^{-1}\left(F-f^{n}e-\Delta
t\tA \left(\frac{1}{\varepsilon}
G_P(F)+L(F)\right)\right)\] which gives \[
f^{n+1}=f^{n}\left(1-w^{T}A^{-1}e\right)+\Delta
t\left(\tw^{T}-w^{T}A^{-1}\tA\right)\left(\frac{1}{\varepsilon}
G_P(F)+L(F)\right)+{w}^{T}A^{-1}F.\]
Using now the assumption (\ref{eq:ai1}) we get
\be
f^{n+1}=f^{n}\left(1-w^{T}A^{-1}e\right)+w^{T}A^{-1}F,
\ee
which permits to pass to the limit $\varepsilon\to 0$ in (\ref{eq:GIMEXbv})-(\ref{eq:GIMEX1bv}).\\
$\Box$

Note that condition (\ref{eq:ai1}) is automatically satisfied if the IMEX scheme is GSA. In this case it is easy to show that we have
\begin{theorem}
If the penalized IMEX scheme is of type A and GSA then
\begin{equation}
\lim_{\varepsilon\to 0} f^{n+1}=M[f^{n+1}].
\end{equation}
\end{theorem}

We consider now the case of penalized IMEX schemes of type CK. We can state an analogous result of Theorem \ref{th:apCK} for standard IMEX schemes of type CK. Some relevant differences will be pointed out at the end of the proof. We have the following

\begin{theorem}
\label{th:papCK} If the penalized IMEX scheme is of type CK and GSA then for consistent initial data in the limit $\varepsilon\rightarrow 0$ scheme
(\ref{eq:GIMEXbv})-(\ref{eq:GIMEX1bv}) becomes the explicit RK scheme characterized by
($\tA, \tw, \tc$) applied to the limit
Euler system (\ref{eq:Euler}).
\end{theorem}

\proof
Using the same notations as in Theorem \ref{th:apCK} the penalized IMEX schemes (\ref{eq:GIMEXbv})-(\ref{eq:GIMEX1bv}) can be written as
\begin{eqnarray}
\nonumber
F^{(1)}&=&f^{n}\\
\hat{F} &=& \displaystyle f^{n}\hat{e}+\Delta t \ta \left(\frac{1}{\varepsilon}
G_P(F^{(1)})+L(F^{(1)})\right)+\Delta t
\hat{\tA}\left(\frac{1}{\varepsilon}
G_P(\hat{F})+L(\hat{F})\right)
\label{eq:theo3ab}
\\
\nonumber
&+&\frac{\Delta
t}{\varepsilon} a Q_P(F^{(1)})+
\frac{\Delta t}{\varepsilon}
\hat{A}Q_P(\hat{F}),\\
\nonumber
f^{n+1}&=&f^n +\Delta t\tw_1 \left(\frac{1}{\varepsilon}
G_P(F^{(1)})+L(F^{(1)})\right)+\Delta t\hat{\tw}^T \left(\frac{1}{\varepsilon}
G_P(\hat{F})+L(\hat{F})\right)\\[-.25cm]
\label{eq:theo3bb}
\\[-.25cm]
\nonumber
&+&w_1\frac{\Delta t}{\varepsilon} Q_P(F^{(1)})+\hat{w}^T\frac{\Delta t}{\varepsilon} Q_P(\hat{F}),
\end{eqnarray}
whereas the corresponding moment scheme is unchanged and reads
\begin{eqnarray}
\nonumber
\langle \phi F^{(1)}\rangle &=& \langle \phi f^{n}\rangle\\[-.3cm]
\label{eq:GIMEXaab} \\[-.3cm]\nonumber\displaystyle   \langle \phi \hat{F}\rangle &=& \langle \phi f^{n}\hat{e}\rangle+\Delta t \ta \langle\phi L(F^{(1)})\rangle
+\Delta t
\hat{\tA}\langle\phi L(\hat{F})\rangle\\
\displaystyle   \langle \phi f^{n+1}\rangle &=& \langle\phi f^{n}\rangle+\Delta t\tw_1\langle\phi L(F^{(1)})\rangle+\Delta t\hat{\tw}^T \langle\phi L(\hat{F})\rangle.
\label{eq:GIMEXbbb}
\end{eqnarray}

Solving now,
the second equation in (\ref{eq:theo3ab}) for $Q_P(\hat{F})$ we get \begin{eqnarray}
\nonumber
\Delta t Q_P(\hat{F})&=&\varepsilon
\hat{A}^{-1}\left[\hat{F}-f^{n}\hat{e}-\Delta t
\hat{\tA}\left(\frac{1}{\varepsilon}
G_P(\hat{F})+L(\hat{F})\right)-\Delta t \ta
\left(\frac{1}{\varepsilon}
G_P(f^n)+L(f^n)\right)\right]\\[-.25cm]
\label{eq:ll1b}
\\[-.25cm]
\nonumber
&-&\Delta t
\hat{A}^{-1}a Q_P(f^{n}).
\end{eqnarray}
As $\varepsilon\to 0$ we obtain
\be
Q_P(\hat{F})= -\hat{A}^{-1}\left(\hat{\tA}G_P(\hat{F})+\ta G_P(f^n)\right)-
\hat{A}^{-1}a Q_P(f^{n}).
\label{eq:limitb}
\ee
Now since $f^n$ is consistent in the limit
$\varepsilon\rightarrow 0$ we have $f^n=M[f^n]$, $Q_P(f^n)=0$, $G_P(f^n)=0$ and (\ref{eq:limitb}) reduces to
\be
Q_P(\hat{F})=-\hat{A}^{-1}\hat{\tA}G_P(\hat{F})
\label{eq:limitb2}
\ee
which thanks to the fact that $\hat{A}^{-1}\hat{\tA}$ is lower triangular with zero diagonal elements implies $\hat{F}=M[\hat{F}]$.
Moreover, since the scheme is GSA, we also have $f^{n+1}=M[f^{n+1}]$ and thus at the next time step the initial value remains consistent and the moments system (\ref{eq:GIMEXaab})-(\ref{eq:GIMEXbbb}) corresponds to the explicit Runge-Kutta methods for the Euler equations.\\
$\Box$

Note that even for ARS type IMEX schemes with $a=0$ the diagonal elements of the right hand side of equation (\ref{eq:limitb}) do not vanish and thus asymptotic preservation is not achieved unless the initial data is consistent or the following additional condition holds true
\be
-\hat{A}^{-1}\ta = 0.
\ee
In this latter case, ARS type schemes are asymptotically accurate if $\tw_1=0$.

Let us consider now the behavior of the penalized IMEX schemes of type CK in case the GSA and the consistent initial data assumptions are removed.
The numerical solution takes the form
\begin{eqnarray*}
f^{n+1}&=&f^n\left(1-\hat{w}^T
\hat{A}^{-1}\hat{e}\right)+\Delta t\left(\tw_1-\hat{w}^T
\hat{A}^{-1}\hat{\ta}\right) \left(\frac{1}{\varepsilon}
G_P(f^{n})+L(f^{n})\right)+\hat{w}^T
\hat{A}^{-1}\hat{F}\\&+&\Delta t\left(\hat{\tw}^T-\hat{w}^T
\hat{A}^{-1}\hat{\tA}\right) \left(\frac{1}{\varepsilon}
G_P(\hat{F})+L(\hat{F})\right)+\frac{\Delta t}{\varepsilon}\left(w_1-\hat{w}^T\hat{A}^{-1}a\right) Q_P(f^n),
\end{eqnarray*}
and thus we need the further assumptions
\be
\tw_1=\hat{w}^T
\hat{A}^{-1}\hat{\ta},\qquad \hat{\tw}^T=\hat{w}^T
\hat{A}^{-1}\hat{\tA},\qquad w_1=\hat{w}^T\hat{A}^{-1}a,
\ee
in order to pass to the limit $\varepsilon\to 0$. Clearly the assumptions above are satisfied by GSA schemes.

We end the paragraph by emphasizing the following result on the asymptotic behavior of the numerical solution
\begin{theorem}
If the penalized IMEX scheme is of type CK and GSA then
\be
\lim_{\varepsilon\to 0} f^{n+1}=M[f^{n+1}],
\ee
if one of the following conditions is satisfied
\begin{itemize}
\item[\rm (a)] the initial data is consistent;
\item[\rm (b)]
\be
\hat{e}_{\nu}^T\hat{A}^{-1}\hat{\tA}=0,\qquad \hat{e}_{\nu}^T\hat{A}^{-1}\ta=0,\qquad \hat{e}_{\nu}^T\hat{A}^{-1}a=0,
\label{eq:apars}
\ee
where $\hat{e}_{\nu}=(0,\ldots,0,1)^T\in\R^{\nu-1}$;
\item[\rm (c)]
\be
\hat{A}^{-1}\ta=0,\qquad \hat{A}^{-1}a=0.
\label{eq:apars2}
\ee
\end{itemize}
\end{theorem}

\subsection{Properties of penalized IMEX-BGK schemes for the Boltzmann equation}
In this last paragraph, we restrict to the particular case where the penalization term is given by a BGK relaxation operator. Exactly as for standard IMEX schemes a fundamental property of equations
(\ref{eq:GIMEXbv})-(\ref{eq:GIMEX1bv}) with $Q_P(f)=\mu(M[f]-f)$ is that they can be solved explicitly.

In fact, since the implicit scheme is a DIRK method the stage values take the form
\begin{eqnarray}
   F^{(i)} = \displaystyle f^{n}+\Delta t \sum_{j=1}^{i-1} \ta_{ij} \left(L(F^{(j)})+\frac{1}{\varepsilon}
G_P(F^{j})\right)+\Delta t\sum_{j=1}^{i} a_{ij}\frac{\mu}{\varepsilon}(M[F^{(j)}]-F^{(j)})\label{eq:GIMEXdirkb}
\end{eqnarray}
where the only implicit term is the diagonal factor $M[F^{(i)}]-F^{(i)}$ in which $M[F^{(i)}]$ depends only on the moments $\langle\phi F^{(i)}\rangle$. If we now integrate equation (\ref{eq:GIMEXdirkb}) agains the collision invariants thanks to the conservations (\ref{eq:QC}) we obtain again the moment scheme (\ref{eq:GIMEXdirkm}). Thus $\langle \phi F^{(i)}\rangle$, and so $M[F^{(i)}]$, can be explicitly evaluated and system (\ref{eq:GIMEXdirkb}) is explicitly solvable. This property was at the basis of the penalized IMEX schemes recently proposed in~\cite{Filbet}.

Next we focus on the monotonicity properties of the penalized IMEX schemes. It appears clear that in space non homogeneous situations due to the penalization approach it is extremely difficult to achieve positivity of the numerical solution. For this reason we restrict to the space homogeneous case, which in this case however involves the full penalized IMEX method. As usual, in combination with time splitting strategies, these results permit to recover
positivity even in the non homogeneous case if required. This is of paramount importance, for instance, in the case of Monte Carlo methods~\cite{PRHawaii}.


In the homogeneous case the method takes the form
\begin{eqnarray}
\label{eq:s1b}
F=f^{n}e+\frac{\mu\Delta t}{\varepsilon}\tA \left(\frac{P(F)}{\mu}- M[f^n] e\right)+\frac{\mu\Delta t}{\varepsilon}A\left(M[f^n]e-F\right),\\
\label{eq:s2b}
f^{n+1}=f^n+\tw^T\frac{\mu\Delta t}{\varepsilon}\left(\frac{P(F)}{\mu}- M[f^n] e\right)+w^T\frac{\mu\Delta t}{\varepsilon}\left(M[f^n]e-F\right),
\end{eqnarray}
where again the local equilibrium $M[F]=M[f^n]e$ is independent of time since $\langle\phi f^{n+1}\rangle=\langle\phi f^n\rangle$. In (\ref{eq:s1b})-(\ref{eq:s2b}) we used the fact that $G_P(F)=P(F)-\mu M[f^n]e$ with $P(F)=Q_B(F)+\mu F \geq 0$ and $\langle \phi P(F)\rangle=\langle \phi \mu F\rangle$.

Let us
define $z=\mu\Delta t/\varepsilon$ and solve for $F$, we get \begin{eqnarray}
F&=&(I+zA)^{-1}\left(f^{n}e+z(A-\tA)M[f^n]e+z\tA \frac{P(F)}{\mu}\right)\\
f^{n+1}&=&(1-w^{T}A^{-1}e)f^{n}+w^{T}A^{-1}F+z(\tw^T-w^TA^{-1}\tA)\left(\frac{P(F)}{\mu}- M[f^n]e\right),
\end{eqnarray}
where we used the identity
\be
z(M[f^n]e-F)=A^{-1}\left[F-f^{n}e-z\tA \left(\frac{P(F)}{\mu}- M[f^n]e\right)\right].\label{eq:simpb}\ee
We can thus state the following
\begin{proposition}
Sufficient conditions to guarantee that $f^{n+1}\geq 0$ when $f^n\geq 0$ in (\ref{eq:s1b})-(\ref{eq:s2b}) are that
\begin{eqnarray}
\label{eq:posnb}
\left(I+zA\right)^{-1}e\geq 0,&\qquad \left(I+zA\right)^{-1}(A-\tA)e\geq 0,\qquad& \left(I+zA\right)^{-1}\tA\geq 0\\
1-w^{T}A^{-1}e \geq 0,&\qquad w^{T}A^{-1}\geq 0,\qquad& \tw^T=w^TA^{-1}\tA.
\label{eq:posn1b}
\end{eqnarray}
\end{proposition}
Conditions (\ref{eq:posnb})-(\ref{eq:posn1b}) applies to a general
penalized IMEX method and must be interpreted component by
component. Clearly for GSA schemes conditions reduce to
(\ref{eq:posnb}). Since (\ref{eq:posnb})-(\ref{eq:posn1b}) are
originated from the penalizing strategy in the nonlinear case they
differ from the classical conditions for absolute monotonicity of
additive Runge-Kutta methods~\cite{Hi1}. The range of values of $z$
for which (\ref{eq:posnb})-(\ref{eq:posn1b}) hold true defines the
\emph{absolute monotonicity region} of the penalized IMEX method.
Note that here we are interested in the behavior of the schemes
outside the stability region of explicit methods, thus monotonicity
is important for $z\geq 1$.

Finally, in order to characterize the stability property of the penalized IMEX schemes we consider the linearized problem
\be
\partial_t f = \frac{1}{\epsilon}G_L(f)+\frac{\mu}{\varepsilon}(M[f]-f),
\ee
where $G_L(f)=(\lambda-\mu)(M[f]-f)$. Here $\mu$ is an estimate of the unknown relaxation rate $\lambda>0$. We will further denote by $\alpha=\lambda/\mu>0$ the penalization parameter of the method. Clearly $\alpha=1$ corresponds to the optimal choice for which the whole analysis reduces to the one performed in Section \ref{par:BGK} for the BGK case.

The penalized IMEX schemes take the form
\begin{eqnarray}
\label{eq:s1bl}
F=f^{n}e+z\tA (\alpha-1)\left(F-M[f^n] e\right)+zA\left(M[f^n]e-F\right),\\
\label{eq:s2bl}
f^{n+1}=f^n+\tw^Tz(\alpha-1)\left(F-M[f^n] e\right)+w^Tz\left(M[f^n]e-F\right).
\end{eqnarray}
Now since
\be
(F-M[f^n]e)=(I+z(A-(\alpha-1)\tA))^{-1}(f^{n}-M[f^n])e.\label{eq:simppl}\ee
the numerical solution becomes
\be
f^{n+1}=R(\alpha,z)f^{n}+(1-R(\alpha,z))M[f^n],
\label{eq:dirksl}
\ee
where \be
R(\alpha,z)=1-z(w^{T}-(\alpha-1)\tw^T)(I+z(A-(\alpha-1)\tA))^{-1}e
\label{eq:fstab}
\ee
plays the role of the stability function of the penalized IMEX method. In this case $R(\alpha,z)$ is an approximation of $\exp(-\alpha z)$ and we characterize the absolute stability region by the condition $|R(\alpha,z)|\leq 1$. Note that $R(1,z)$ corresponds to the  stability function of the implicit method.


In particular, similar to the BGK case (\ref{eq:Rrec}), we
recursively obtain \be
f^{n+1}=R(\alpha,z)^{n+1}f^{n}+(1-R(\alpha,z)^{n+1})M[f_0]. \ee This
permits to consider a weaker notion of AP property when
$|R(\alpha,\infty)|<1$ where \be
R(\alpha,\infty)=\lim_{z\to\infty}R(\alpha,z). \label{eq:fstabl} \ee
In this case intact $f^{n+1}\to M[f_0]$ as $n\to\infty$ at a rate
characterized by the damping constant $|R(\infty,\alpha)|<1$. Of
course for an L-stable implicit method we have $R(1,\infty)=0$. From
(\ref{eq:fstabl}) an AP method corresponds to the requirement
$R(\alpha,\infty)=0$, $\forall\, \alpha>0$. A condition which
essentially generalizes the notion of L-stability to penalized IMEX
schemes.

The sensitivity of $R(\alpha,\infty)$ to values of the penalization
parameter $\alpha\neq 1$ gives a measure of the robustness of the
method with respect to the estimate $\mu$ of $\lambda$ in terms of
capturing the correct asymptotic behavior. Thus we introduce the
following
\begin{definition}
The values of $\alpha>0$ for which $|R(\alpha,\infty)|<1$ characterizes the \emph{weak AP range} of the penalized IMEX method.
\label{def:wap}
\end{definition}

\begin{remark}
Note that, to achieve better penalization properties, the value of
$\mu$ can be taken space and time dependent in (\ref{eq:GIMEXdirkb})
to minimize the value of $G_P(F^{(j)})$. Typically, for a given $f$,
one should choose $\mu>0$ such that it locally minimizes $\|P(f)-\mu
M[f]\|$ in a suitable norm $\|\cdot\|$. Here we do not explore
further this aspect. However, the general asymptotic preserving
properties and considerations we have done remain valid also in this
case.
\end{remark}

%

\section{Numerical considerations and tests}

We report in Table \ref{tab:1} a summary of the properties of some
penalized IMEX schemes satisfying the GSA property developed in the
literature. We use the notation {NAME$(\nu_E,\nu_I,p)$} where
$\nu_E,\nu_I$ are, respectively, the number of function evaluations of the
explicit and the implicit methods and $p$ is the combined order of
the IMEX scheme. The field NAME of the schemes is composed by the
initials of the authors and the scheme type. We consider the
following methods: ARS$(1,1,1)$ corresponding to forward-backward
Euler, ARS$(2,2,2)$ and ARS$(4,4,3)$ from Sections 2.6 and 2.8
in~\cite{Ascher}, DP-ARS$(2,2,2)$ same as ARS$(2,2,2)$ but with
$\gamma=(2+\sqrt{2})/2$, JF-CK$(2,3,2)$ from (2.8) Section 2
in~\cite{Filbet}, BPR-CK$(3,5,3)$ from the Appendix in~\cite{BPR}.
Schemes DP1-A$(1,2,1)$, DP2-A$_2$(1,2,1), DP-ARS$(1,2,1)$ and DP-A$(2,4,2)$
are reported in a separate Appendix at the end of the manuscript.
With the acronym AA we denote the asymptotic accurate property, with
AA-c the same property but for a consistent initial data. We also
report the region of absolute monotonicity (AM) in the nonlinear
case accordingly to definition (\ref{eq:posnb})-(\ref{eq:posn1b}),
the asymptotic behavior of the stability function (\ref{eq:fstabl})
and the weak AP range of the schemes obtained from Definition
\ref{def:wap}. For schemes marked with an '*' the DIRK part has
order $p+1$.

\begin{table}[t]
\caption{Summary of the properties of penalized IMEX GSA schemes}
\begin{center}
{\small
\begin{tabular}{c|c|c|c|c|c}
\hline
&&&&&\\[-.2cm]
Name & AA & AA-c & AM & $R(\alpha,\infty)$& $|R(\alpha,\infty)|<1$\\[+.1cm]
\hline
&&&&&\\[-.2cm]
ARS$(1,1,1)$ & no & yes & $z\in[0,\infty)$&  $\alpha-1$ & $\alpha\in(0,2)$\\[+.1cm]
\hline
&&&&&\\[-.2cm]
DP-ARS$^*(1,2,1)$ & yes & yes & $z=0$&  $0$ & $\alpha\in(0,\infty)$\\[+.1cm]
\hline
&&&&&\\[-.2cm]
DP-A$^*(1,2,1)$ & yes & yes & $z\in[0,\infty)$&  $0$ & $\alpha\in(0,\infty)$\\[+.1cm]
\hline
&&&&&\\[-.2cm]
ARS$(2,2,2)$ & no & yes & $z=0$&  $\frac{(1-\alpha)(\alpha + 4\gamma - 2\alpha\gamma - 2\gamma^2)}{2\gamma^2}$ & $\alpha\in(0.874,1.117)$\\[+.1cm]
\hline
&&&&&\\[-.2cm]
DP-ARS$(2,2,2)$ & no & yes & $z\in[0,\infty)$ & $\frac{(1-\alpha)(\alpha + 4\gamma - 2\alpha\gamma - 2\gamma^2)}{2\gamma^2}$ & $\alpha\in(0,2.288)$\\[+.1cm]
\hline
&&&&&\\[-.2cm]
JF-CK$(2,3,2)$ & no & yes & $z\in[0,2]$ & {\footnotesize $2\alpha^2-4\alpha+1$} & $\alpha\in(0,2)$\\[+.1cm]
\hline
&&&&&\\[-.2cm]
DP1-A$^*(2,4,2)$ & yes & yes & $z=0$ & $0$ & $\alpha\in(0,\infty)$\\[+.1cm]
\hline
&&&&&\\[-.2cm]
DP2-A$_2(2,4,2)$ & yes & yes & $z\in[1,\infty)$ & $0$ & $\alpha\in(0,\infty)$\\[+.1cm]
\hline
&&&&&\\[-.2cm]
ARS$(4,4,3)$ & no & yes & $z=0$ & $\frac{(1-\alpha)(7\alpha^3 - 78\alpha^2 + 138\alpha - 38)}{18}$ & $\alpha\in(0.13475,2)$\\[+.1cm]
\hline
&&&&&\\[-.2cm]
BPR-CK$(3,5,3)$ & no & yes & $z=0$ & $\frac{4(2\alpha^3 - 5\alpha^2 + 2\alpha)}{3} + 1$&{{$\alpha\in(0,\frac{1-\sqrt{3}}{2})\cup(\frac12,\frac{1+\sqrt{3}}{2})\cup(\frac32,2)$}}\\[+.1cm]
\hline
\end{tabular}
}
\end{center}
\label{tab:1}
\end{table}%

We emphasize that the computational cost of penalized IMEX schemes is characterized by the number of
stages of the explicit method since, by construction, the implicit part is applied to the easy
invertible term used for penalization. This is an important aspect for
the practical construction of the methods.

Next we numerically measure the accuracy of second and third order methods in solving the full
Boltzmann equation (\ref{eq:1})-(\ref{eq:Q}) on a periodic smooth
solution. We emphasize that to our knowledge these are the first full third order implementations for the Boltzmann equation.

The computation is performed on $(x, v) \in [0, 1] \times
[-v_{\max}, v_{\max}]^2$, with $v_{\max} = 8$. We use a $3$rd order WENO
scheme for the space discretization \cite{Shu} and a fast spectral
method for solving the collision integral \cite{MP}. We take $N_v = 32$
grid points in each velocity direction which provides enough accuracy in $v$ to measure the convergence rates in time. The time step
is fixed equal to $\Delta t = {\Delta x}/({2v_{\max}})$. The initial data is \be
\varrho_{0}(x)=\frac{2+\sin(2\pi x)}{3}, \ u_{0}(x)=\frac{\cos(2 \pi
x)}{5}, \ T_{0}(x)=\frac{3+\cos(2 \pi x)}{4}.\ee We report the $L_1$
norm of the error for the density for different values of the
Knudsen number, i.e. $\varepsilon=10^{-1}$, $\varepsilon=10^{-3}$
and $\varepsilon=10^{-6}$. We consider both equilibrium initial
data $f_0(x,v)=M[f_0]$ and non equilibrium initial
data \be f_0(x,v)=\frac{\varrho_{0}(x)}{(2\pi
T_{0}(x))^{1/2}}\frac{1}{2}\left(\exp^{-\frac{|v-u_0(x)|^{2}}{2T_0(x)}}+\exp^{-\frac{|v+3u_0(x)|^{2}}{2T_0(x)}}\right).\ee
To measure the convergence rate we repeat the computation for an
increasing number of grid points in space, $N_x=128, 256, 512,
1024$. 

\begin{figure}[h!]
\begin{center}
\includegraphics[scale=0.38]{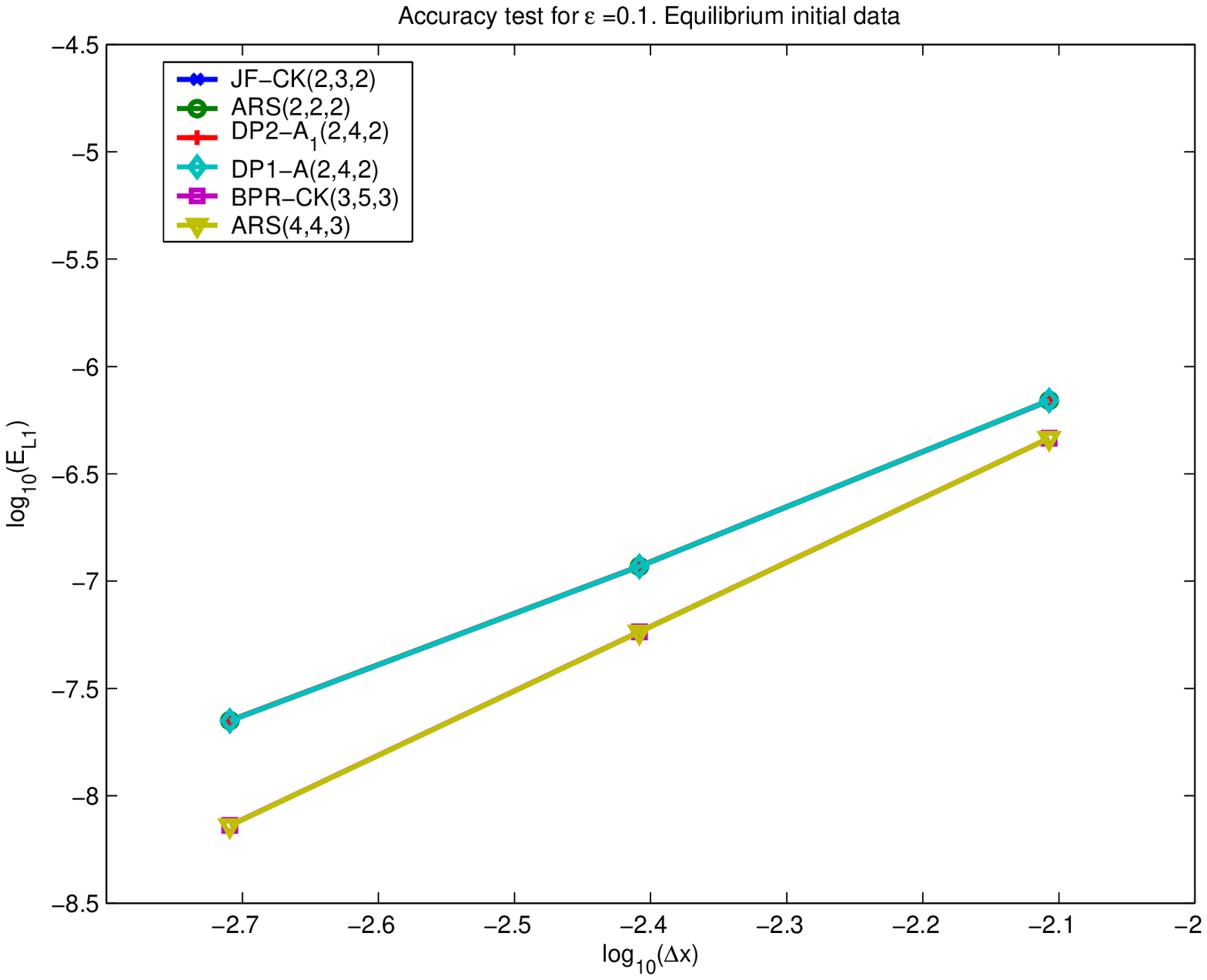}
\includegraphics[scale=0.38]{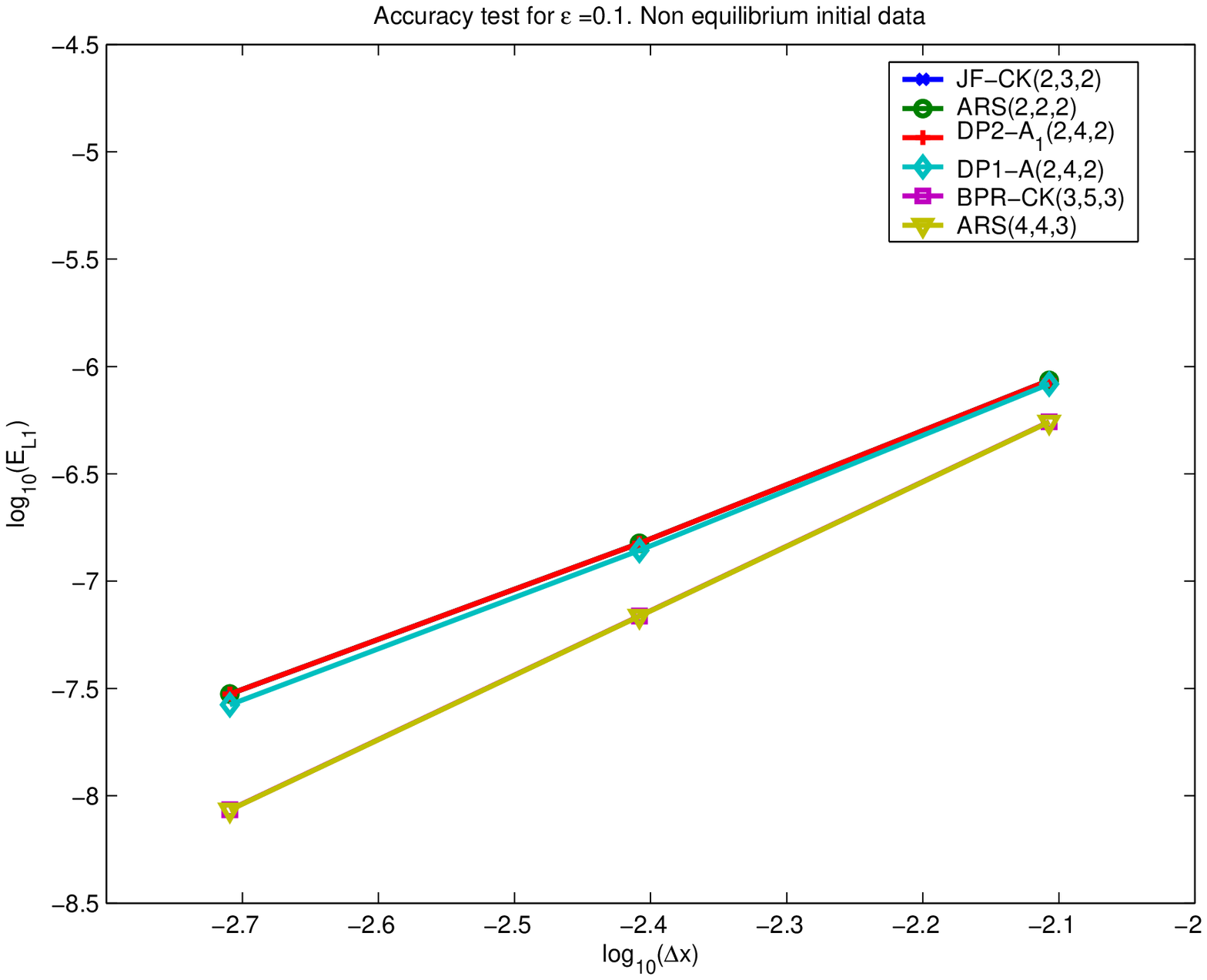}\\
\includegraphics[scale=0.38]{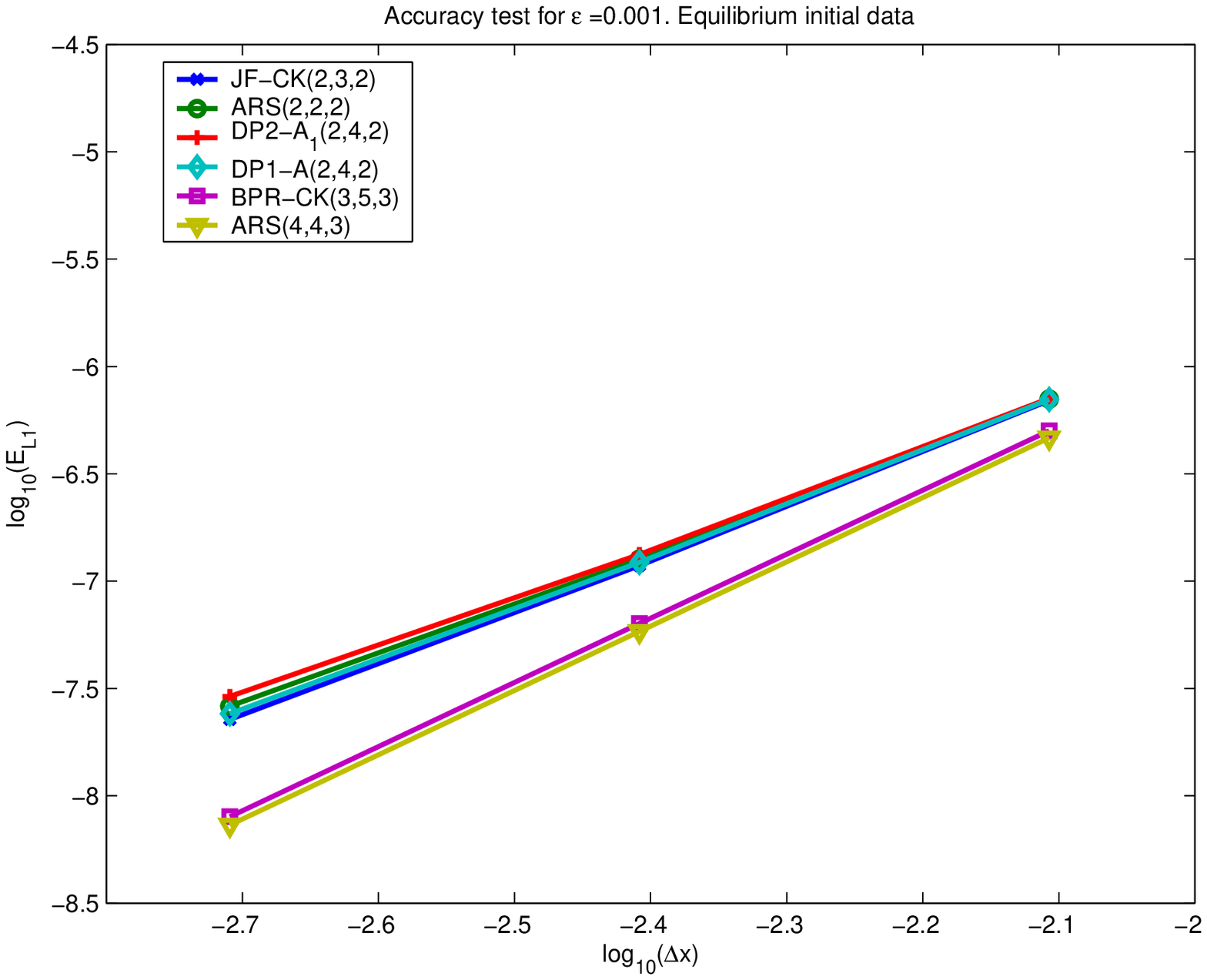}
\includegraphics[scale=0.38]{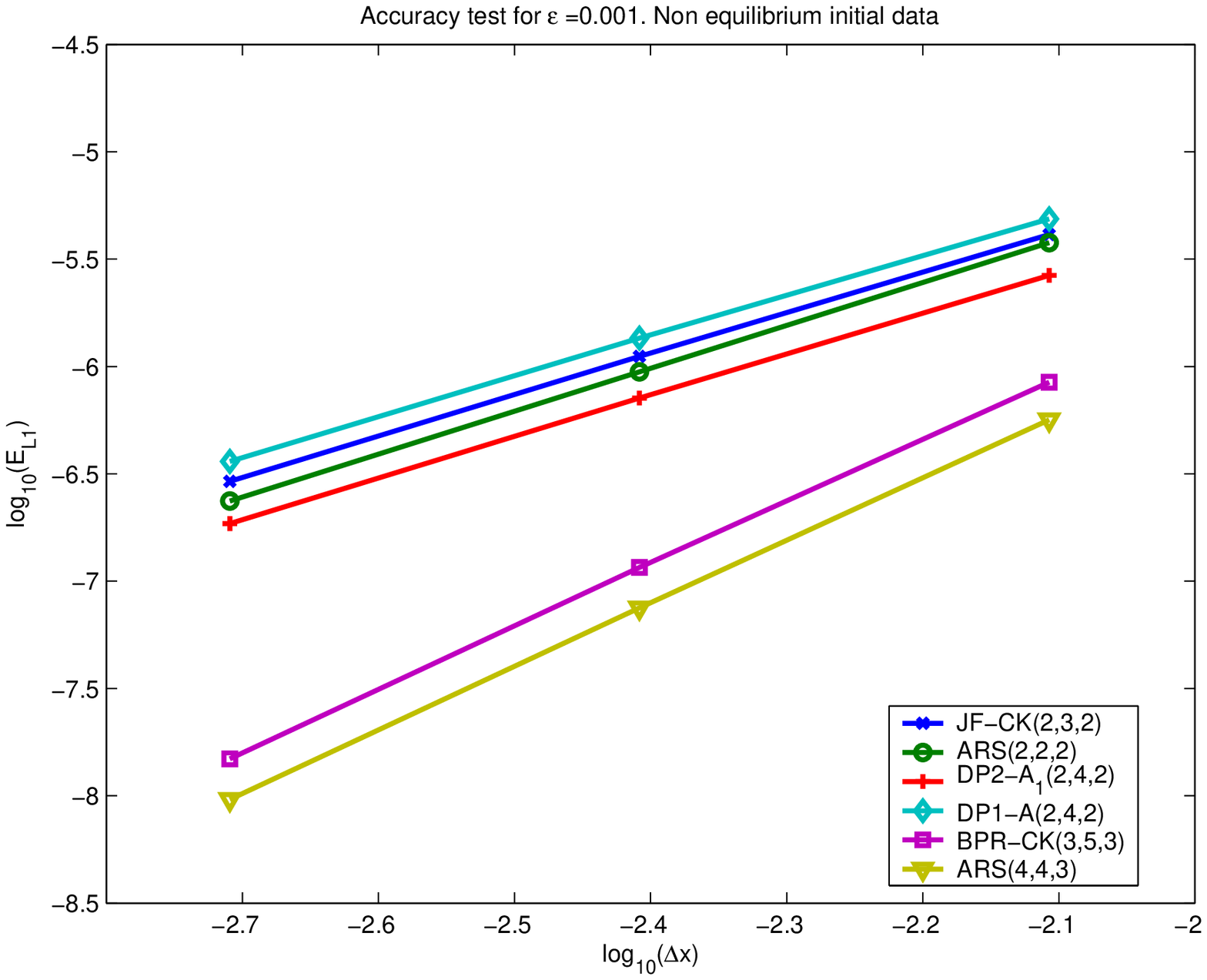}\\
\includegraphics[scale=0.38]{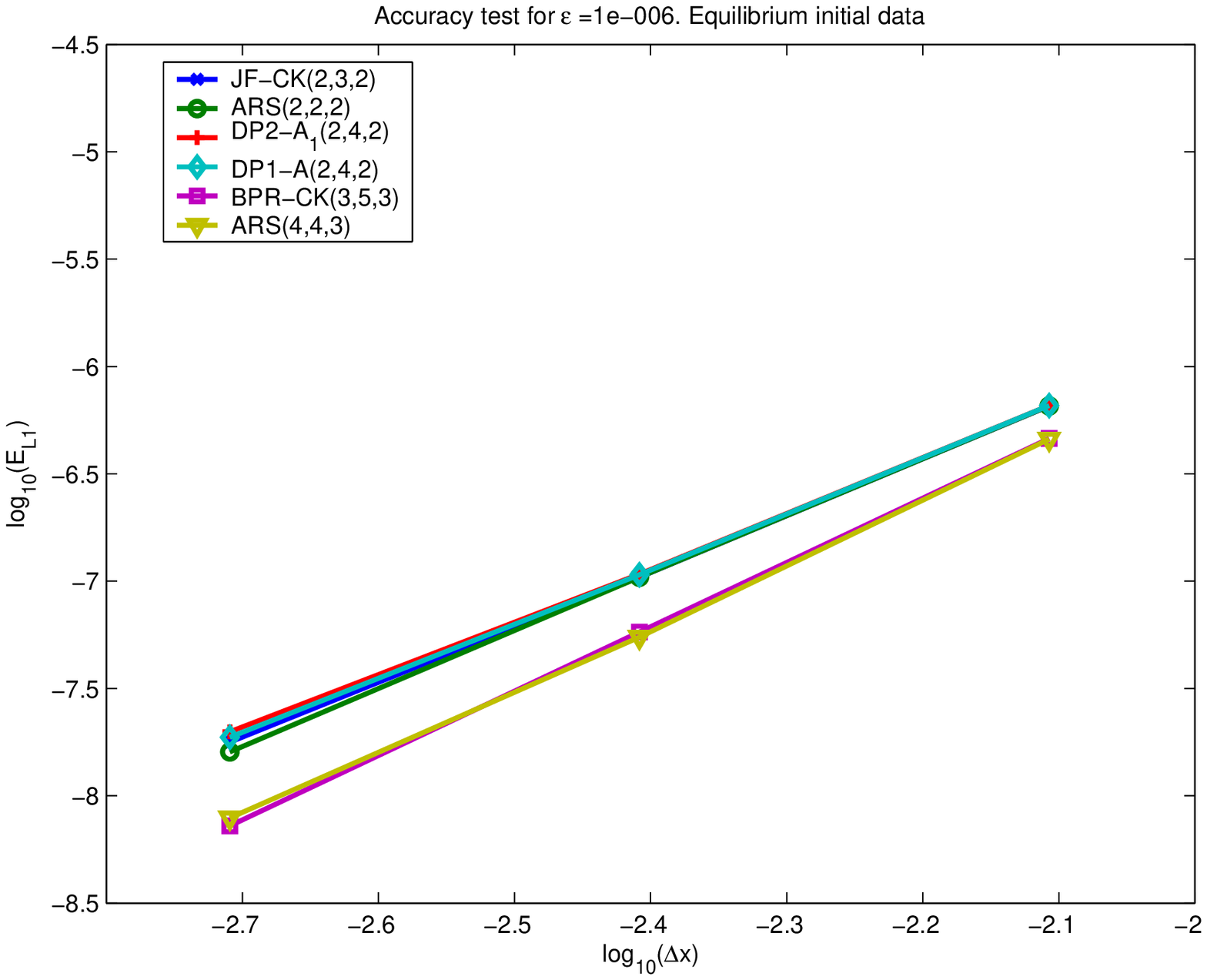}
\includegraphics[scale=0.38]{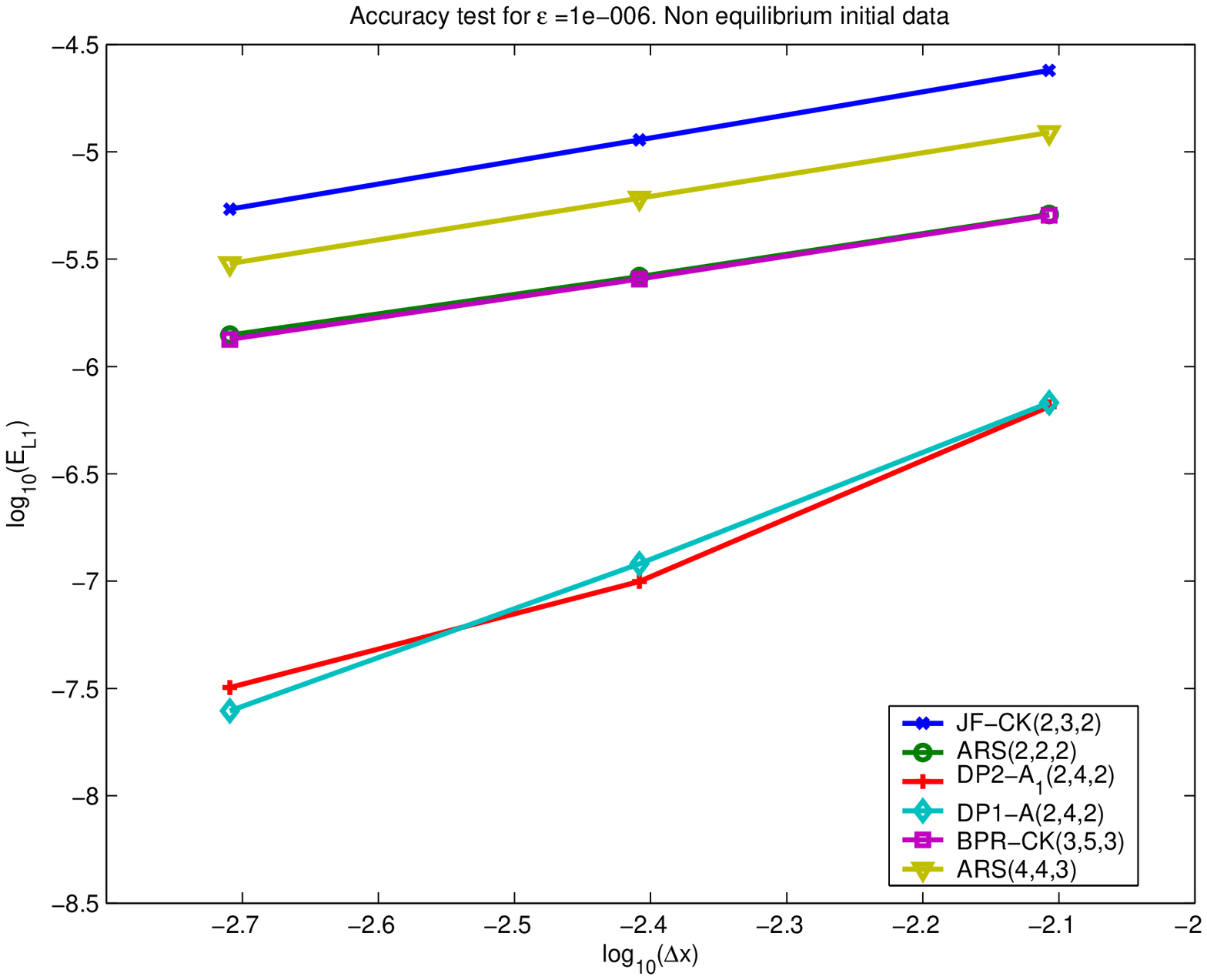}
\caption{$L_1$ error for the density $\varrho$ for different second
and third order IMEX schemes. Left column equilibrium initial data,
right column non equilibrium initial data. Top
$\varepsilon=10^{-1}$, center $\varepsilon=10^{-3}$, bottom
$\varepsilon=10^{-6}$.}\label{fig:conv}
\end{center}
\end{figure}

\begin{figure}[h!]
\begin{center}
\includegraphics[scale=0.38]{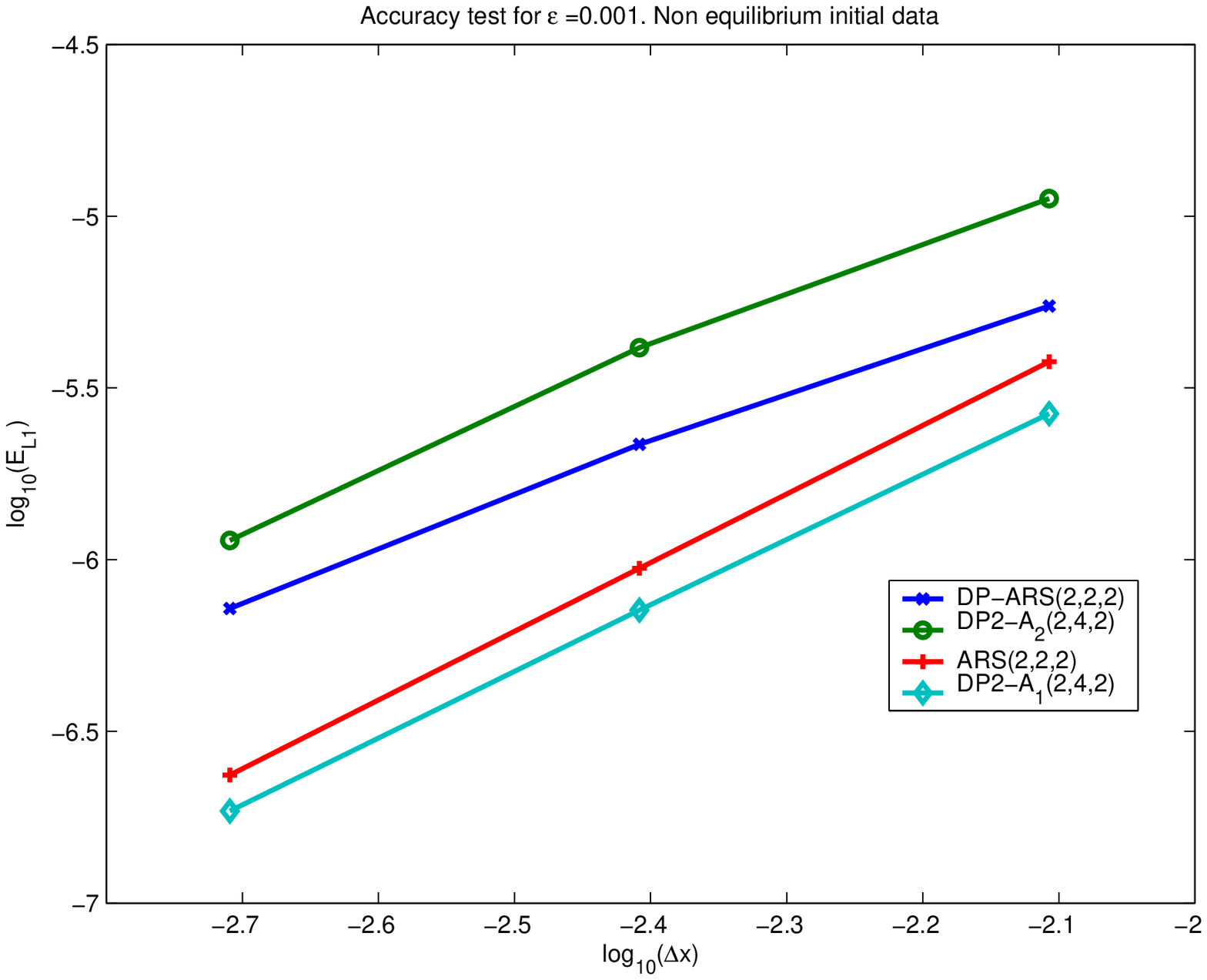}
\caption{$L_1$ error for the density $\varrho$ for second order monotone schemes for different values of the diagonal entries in the monotone and non monotone regions.
Non equilibrium initial data with
$\varepsilon=10^{-3}$.}\label{fig:conv2}
\end{center}
\end{figure}

Figure \ref{fig:conv} shows the results for the second and third order penalized IMEX schemes. Here in scheme DP2-A(2,4,2) we choose to maximize accuracy taking $\gamma=1/3$ outside the monotonicity region. As
expected, all the schemes exhibit the prescribed order of
convergence for equilibrium initial data while 
degradation of accuracy to first order is observed for schemes not asymptotically accurate. Schemes DP1-A$(2,4,2)$
and DP2-A$(2,4,2)$ preserves second order accuracy uniformly in all regimes independently of the initial data. Finally, we report in figure \ref{fig:conv2} the
comparisons between the convergence rate of the second order monotone schemes where the diagonal value is taken inside and outside the monotonicity region. Note how the large diagonal entry needed for monotonicity reduces the accuracy of the schemes for non equilibrium initial data in the intermediate regimes.

\section{Conclusions and further developments}

We have studied the extension of implicit-explicit (IMEX) Runge-Kutta
schemes to the numerical solution of nonlinear kinetic equations in
stiff regimes. We derived sufficient conditions
for an IMEX scheme to be asymptotic preserving and asymptotically
accurate. These notions are of paramount importance in the practical computation of kinetic equations where regions with a large variation of the computational space and time scales are present.

In the first part we focused on standard IMEX methods and studied their asymptotic behavior and monotonicity property. In the second part, following the approach introduced in \cite{Filbet}, we tackle the challenging case of the efficient computation of the Boltzmann collision term in stiff regimes. A penalization technique permits to derive new schemes, here referred to as penalized IMEX schemes, which under additional assumptions keep the asymptotic properties of the standard IMEX methods by avoiding the costly inversion of the collision term. Monotonicity is studied for the case where the penalization term is represented by a relaxation operator.

For the future extensive numerical testing and applications of the schemes to other collisional
kinetic equations, like the Landau-Fokker-Planck equation, is
scheduled.

\vskip 1cm \textbf{Acknowledgement}. G. Dimarco was supported by the
French ANR project BOOST. The authors would like to thank Prof.
S. Jin and Prof. F. Filbet for stimulating discussions.

\section*{Appendix}
In the sequel we report the Butcher tableu of some penalized IMEX schemes satisfying the AA condition. The guidelines in the design of the methods are the GSA property and the minimization of the number of explicit function evaluations. We also consider schemes for which the explicit part is strongly stability preserving~\cite{GST}. As we will see since the above requirements lead to schemes using several implicit levels, it is natural to satisfy higher order conditions for the DIRK part.

Penalized first order schemes satisfying the AA property are easy to construct, we give here two examples, one absolutely monotone of type A, referred to as DP-A$(1,2,1)$
{\small
\[
    \begin{array}{c|cc}
      0 & 0 & 0 \\
      1   & 1 & 0 \\
      \hline
          & 1 & 0
    \end{array}\qquad
  \begin{array}{c|cc}
      \gamma & \gamma & 0 \\
      1   & 1-\gamma & \gamma \\
      \hline
          & 1-\gamma & \gamma
    \end{array}
    \]}
where $\gamma\geq 1/2$,
and a non monotone variant of type ARS that we refer to as DP-ARS$(1,2,1)$ which satisfies the asymptotic preservation conditions (\ref{eq:apars})
{\small
\[
    \begin{array}{c|ccc}
      0 & 0 & 0 & 0   \\
      \delta   & \delta & 0 & 0 \\
      1   & 1 & 0 & 0\\
      \hline
          & 1 &  0 & 0
    \end{array}\qquad
  \begin{array}{c|ccc}
      0 & 0 & 0 & 0 \\
      \gamma   & 0 & \gamma & 0 \\
      1   & 0 & 1-\gamma & \gamma\\
      \hline
         & 0 & 1-\gamma & \gamma
    \end{array}
    \]
}
where $\delta=\gamma/(1-\gamma)$. For both first order schemes the choice $\gamma=(1\pm\sqrt{2})/2$ yield  second order accuracy for the DIRK part.

Note that even if two levels are required for the penalized schemes to be GSA and achieve AA, they use just one explicit function evaluation.

For GSA type A schemes it is easy to show that we need at least four levels to satisfy the second order accurate requirement and the GSA property. On the other hand four levels do not suffice to have third order accuracy and GSA property.
As examples we report below scheme DP1-A$(2,4,2)$, non monotone, which achieves third order accuracy on the DIRK part
{\small
\[
    \begin{array}{c|ccccc}
      0 & 0 & 0 & 0 & 0   \\
      1/3   & 1/3 & 0 & 0 & 0 \\
      1   & 1 & 0 & 0 & 0 \\
      1   & 1/2 & 0 & 1/2 & 0 \\
      \hline
          & 1/2   &  0 & 1/2 & 0
    \end{array}\qquad
  \begin{array}{c|ccccc}
      1/2 & 1/2 & 0 & 0 & 0   \\
      2/3   & 1/6 & 1/2 & 0 & 0 \\
      1/2   & -1/2 & 1/2 & 1/2 & 0 \\
      1   & 3/2 & -3/2 & 1/2 & 1/2 \\
      \hline
          & 3/2 & -3/2 & 1/2 & 1/2
    \end{array}
    \]
    }
and a second order scheme DP2-A$(2,4,2)$
{\small
\[
    \begin{array}{c|ccccc}
      0 & 0 & 0 & 0 & 0   \\
      0   & 0 & 0 & 0 & 0 \\
      1   & 0 & 1 & 0 & 0 \\
      1   & 0 & 1/2 & 1/2 & 0 \\
      \hline
          & 0   &  1/2 & 1/2 & 0
    \end{array}\qquad
  \begin{array}{c|ccccc}
      \gamma & \gamma & 0 & 0 & 0   \\
      0   & -\gamma & \gamma & 0 & 0 \\
      1   & 0 & 1-\gamma & \gamma & 0 \\
      1   & 0 & 1/2 & 1/2-\gamma & \gamma \\
      \hline
          & 0   &  1/2 & 1/2-\gamma & \gamma
    \end{array}
    \]
}
which is monotone on the interval $z\geq (2\gamma^2-4\gamma+1)^{-1}$ for $\gamma> (2+\sqrt{2})/2$. In particular for $\gamma=2$ monotonicity is obtained in $[1,\infty)$. From the accuracy viewpoint the choice $\gamma=1/3$ originates a DIRK method which satisfies the third order condition $w^T A\, c = 1/6$. We denote by DP2-A$_1(2,4,2)$ the scheme with $\gamma=1/3$ and by DP2-A$_2(2,4,2)$ the scheme with
$\gamma=2$.

\end{document}